\documentclass[preprint,11pt]{article}
\usepackage{amssymb}
\usepackage{mathrsfs}
\usepackage{amsfonts}
\usepackage{graphicx}
\usepackage{colortbl,dcolumn}
\usepackage{amsmath}
\usepackage{psfrag}
\usepackage{booktabs}
\usepackage{array}
\usepackage{cite}
\numberwithin{equation}{section}
\usepackage[top=1.2in, bottom=1.2in, left=0.9in, right=0.9in]{geometry}

\newcommand{\N}{\mathbb{N}}

\renewcommand{\P}{\mathbb{P}}

\newcommand{\dd}{\text{d}}

\newtheorem{assumption}{{Assumption}}[section]
\newtheorem{theorem}{{Theorem}}[section]
\newtheorem{proposition}{{Proposition}}[section]
\newtheorem{lemma}{{Lemma}}[section]

\newtheorem{example}{{Example}}[section]

\begin{document}

%
 \title {An accelerated exponential time integrator for semi-linear stochastic strongly damped wave equation with additive noise
\footnote{R.Q. was supported by  Research Fund for Northeastern University at Qinhuangdao 
(No.XNB201429), Fundamental Research Funds for Central Universities  (No.N130323015),  
Science and Technology Research Founds for Colleges and Universities in Hebei Province
 (No.Z2014040), Natural Science Foundation of Hebei Province, China (No.A2015501102).
 X.W. was supported by Natural Science Foundation of China (No.11301550, No.11171352).}}

\author{
Ruisheng Qi$\,^\text{a}$,  \quad Xiaojie Wang$\,^\text{b}$ \\
\footnotesize $\,^\text{a}$ School of Mathematics and Statistics, Northeastern University at Qinhuangdao, Qinhuangdao, China\\
\footnotesize qirsh@neuq.edu.cn\; and \;qiruisheng123@sohu.com\\
\footnotesize $\,^\text{b}$ School of Mathematics and Statistics, Central South University, Changsha, China\\
\footnotesize x.j.wang7@csu.edu.cn\; and \;x.j.wang7@gmail.com
}
\maketitle
\begin{abstract}\hspace*{\fill}\\
  \normalsize
This paper is concerned with the strong approximation of a semi-linear stochastic  wave equation with strong damping,
driven by additive noise. Based on a spatial discretization performed by a spectral Galerkin method,
we introduce a kind of accelerated exponential time integrator involving linear functionals of the noise.
Under appropriate assumptions, we provide error bounds for the proposed full-discrete scheme. It is shown that
the scheme achieves higher strong order in time direction than the order of temporal regularity of the underlying problem,
which allows for higher convergence rate than usual time-stepping schemes.   For the space-time white noise case
in two or three spatial dimensions, the scheme still exhibits a good convergence performance.
Another striking finding is that, even for the velocity with low regularity the scheme always promises first order strong convergence
in time. Numerical examples are finally reported to confirm our theoretical findings.

\end{abstract}

\section{Introduction}
\label{sec:Introduction}

Great attention has been devoted in the last decades to numerical approximations of evolutionary stochastic partial
differential equations (SPDEs) (see, e.g. \cite{wang2014higher,Walsh2006numerical,Quer2006space,kovacs2012BITweak,kovacs2013BITweak,Jiang2015stochastic,Cao2007spectral,cohen2015fully} and references therein).
In the present work, we concentrate on a class of semi-linear SPDEs of second order with damping, 
described by
\begin{align}\label{eq:stocha-str-dam-wave-equation}
  \left\{\begin{array}{ll}
\dd u_t =  \alpha L u_t \, \dd t + L u \, \dd t + F(u) \, \dd t +  \dd W(t), & \text{ in }\; \mathcal{D}\times (0,\; T],\\
 u(\cdot,0)=u_0,\;u_t(\cdot,0)=v_0,& \text{ in } \; \mathcal{D},\\
  u=0,& \text{ on } \; \partial \mathcal{D}\times (0,\; T],
  \end{array}\right.
\end{align}
where $\mathcal{D}\subset\mathbb{R}^d$, $d=1,2,3$, is a bounded open domain with a  boundary $\partial\mathcal{D}$,
and where $L := \sum_{i, j = 1}^d \tfrac{\partial}{ \partial x_i}  \big(  l_{ij} (x) \tfrac{\partial}{ \partial x_j}  \big),
x \in \mathcal{D} $ is a linear second-order elliptic operator with smooth coefficients $\{l_{ij}\}_{i,j=1}^d$
being uniformly positive definite. Let $\alpha > 0$ be a fixed positive constant and let $\{W(t)\}_{ t\in [0, T] } $ be a  (possibly cylindrical) $Q$-Wiener process defined on a stochastic basis
$(\Omega,\mathcal{F},\mathbb{P},\{\mathcal{F}_t\}_{t \in [0, T]})$ with
respect to the normal filtration $\{\mathcal{F}_t\}_{t \in [0, T] }$. The initial data
$u_0, v_0$ are assumed to be $\mathcal{F}_0$-measurable random variables.

The deterministic counterpart of \eqref{eq:stocha-str-dam-wave-equation}, called strongly damped wave equation (SDWE),
occurs in a wide range of applications such as modeling motion of viscoelastic materials
\cite{fitzgibbon1981strongly,massatt1983limiting,pata2005strongly}.  From both the theoretical and 
numerical point of view, the deterministic problem has been extensively studied (e.g.,\cite{larsson1991finite,thomee2004maximum,kalantarov2009finite}).
However,  the corresponding stochastic strongly damped wave equations are theoretically and numerically far from well-understood from existing literature\cite{da2014stochastic,qi2015error}.
In the classical monograph \cite{da2014stochastic}, a  stochastic strongly damped wave equation with multiplicative noise
was discussed and a unique mild solution was established.  In our recent publication \cite{qi2015error},
we analyzed the  regularity properties of the mild solution to \eqref{eq:stocha-str-dam-wave-equation}
and  examine error estimates of a full discretization, done by the finite element spatial approximation
together with the well-known linear implicit Euler temporal discretization.
It was shown there that, for a certain class of stochastic SDWEs, the convergence rates of the usual full discretization coincide with  the space-time regularity properties of the mild solution (see Theorem 2.4 in \cite{qi2015error}).
In this article, we aim to introduce a so-called accelerated exponential time integrator for the problem \eqref{eq:stocha-str-dam-wave-equation}, which, as we will show later, promises higher convergence order in time than the order of temporal regularity of the underlying problem.

Different from the usual Euler-typr time-stepping schemes using the basic increments of
the driven Wiener process \cite{yan2005galerkin,kruse2014optimal,lord2013stochastic,anton2015full,
cohen2013trigonometric,qi2013full,qi2013weak,wang2015exponential},
the accelerated exponential time integrators rely on suitable linear functionals of the Wiener process
and  usually attain  superconvergence rates in time \cite{jentzen2009overcoming,jentzen2011efficient,wang2015note,wang2014higher}. In 2009,  such scheme was first constructed by Jentzen and Kloeden \cite{jentzen2009overcoming} for semi-linear parabolic SPDEs with additive space-time white noise. The order barrier in the numerical approximation of parabolic SPDEs was overcome and  a strong convergence rate of order $1-\epsilon$  in time was obtained,  for arbitrarily small $\epsilon>0$, unfortunately with seriously restrictive commutativity condition imposed on the nonlinearity (see \cite[Assumption 2.4]{jentzen2009overcoming}). Afterwards, the accelerated schemes were extended to solve a larger class of parabolic SPDEs with more general noise and the error bounds were analyzed under relaxed conditions on the nonlinearity \cite{jentzen2011efficient,wang2015note}.
Furthermore, the accelerated scheme was successfully adapted to solve semilinear stochastic wave equations and the order barrier $\tfrac 1 2$ was went beyond \cite{wang2014higher}.
%

Following the idea of the acceleration technique, we discretize the considered problem
\eqref{eq:stocha-str-dam-wave-equation} in space by a spectral Galerkin method and
in time by an exponential integrator involving linear functionals of the noise.
To analyze the resulting error bounds, we formulate mild assumptions on the
nonlinear mapping $F$ (see Assumption \ref{assum:F-condition}),
to allow for a large class of nonlinear Nemytskii operators. Additionally, we assume the covariance operator $Q \colon L_2(\mathcal{D}) \rightarrow L_2(\mathcal{D}) $ of the Wiener process obeys
\begin{equation}\label{eq:AQ-condition}
\|A^{\frac{\gamma-1}{2}}Q^{\frac{1}{2}}\|_{\mathrm{HS}}<\infty \quad \text{ for some } \;\gamma\in (-1,2],
\end{equation}
which covers both the space-time white noise case and the trace-class noise case. Here $A := - L$ with 
domain  $D(A) = H^2 (\mathcal{D}) \cap H_0^1 (\mathcal{D})$.
Under these assumptions, the continuous problem \eqref{eq:stocha-str-dam-wave-equation}
possesses a unique mild solution with the displacement $u(t), t \in [0, T]$ taking values in
$L^p (\Omega; \dot{H}^{1+\min\{\gamma, 1\}})$ for $\gamma\in (-1,2]$ and the velocity $v (t), t \in [0, T]$
taking values in  $L^p (\Omega; \dot{H}^{ \gamma }) $ for $\gamma\in [0,2]$.
Our convergence analysis shows that, the convergence rates of the proposed scheme for the displacement $u$ are given by
\begin{align}\label{eq:um-uN}
\|u(t_m)-u^N_m\|_{L^2(\Omega;\dot{H}^0)}
\leq
  C(k^{\min\{1+\gamma,1\}}
  +
  \lambda_{N+1}^{-\frac{1+\min\{\gamma, 1\}}{2}})
  \;
  \text{ for } \;\gamma\in (-1,2],
\end{align}
and for the velocity $v$ in the case $\gamma\in [0,2]$,
\begin{align} \label{eq:vm-vN}
\|v(t_m)-v^N_m\|_{L^2(\Omega;\dot{H}^0)}
\leq
 C(k+\lambda_{N+1}^{-\frac{\gamma}{2}}).
\end{align}
Here $u_m^N$ and $v_m^N$ are the numerical approximations of $u(t_m)$ and $v(t_m)$, respectively.
It should be emphasized that,  although the idea of the method construction comes from \cite{jentzen2009overcoming,jentzen2011efficient,wang2015note,wang2014higher},
the corresponding error analysis (section \ref{sec:proof}) is not easy and
the error analysis
forces us to exploit  a variety of regularity properties (Lemmas \ref{lem:spatial-regualrity-semigroup},
\ref{lemma:result-S4})  of the associated semigroup that are missing in existing works.
As clearly indicated in \eqref{eq:um-uN} and \eqref{eq:vm-vN}, the convergence rates in space
exactly agree with the order of the spatial regularity of the mild solutions. 
However, the convergence rates  in time behave quite differently.
When $\gamma\in(-1,0]$,  the temporal mean-square order
$O(k^{\gamma+1})$ is twice as high as  the temporal H\"older regularity order  of the mild solution 
(see Theorem \ref{thm:main-regularity}). Particularly for the case of the additive space-time white
noise ($Q = I$) in  two and three dimensions when \eqref{eq:AQ-condition} is fulfilled with $\gamma <0 $ and $\gamma <-\tfrac12 $, respectively,  the error estimate \eqref{eq:um-uN} implies that
the introduced method solving the displacement has a good convergence performance.
Also, we point out a surprising fact that,  when $\gamma\in[0,2]$, the numerical solution $u_m^N$ and $v_m^N$
both enjoy a temporal convergence rate of order one that does not depend on the order of the temporal
regularity of the mild solution. This means, even for the velocity $v$ with very low regularity, 
the numerical approximation shows
a first order strong convergence in time. Recall from \cite{qi2015error} that, the finite element spatial
discretization requires \eqref{eq:AQ-condition} to be satisfied with $\gamma \geq 0$ and this excludes
the space-time white noise case in  two or three dimensions.
Furthermore, the convergence rates of numerical approximations obtained in \cite{qi2015error}
coincide exactly with the orders of the regularity. As one can see from earlier discussions,
the newly proposed scheme shows significantly improved convergence rates in time.

Finally, we would like to mention limitations of the presented numerical method.
In order to easily implement the scheme (see subsection 5.1), the eigen-functions of
the dominating linear operator $A$ and of the covariance operator $Q$ of the Wiener process
must coincide and must be known explicitly.


The remainder of this article is organized as follows. In the next section, some preliminaries are collected  and assumptions are made for the noise and the nonlinearity. Further, a concrete example is presented to illustrate the abstract assumptions. In section \ref{sec:main-result} we propose the full-discretization scheme and state the main convergence result, together with some comments on the implementation of the scheme. 
The proof of the main result  is elaborated in section \ref{sec:proof}. 
Finally, numerical experiments are performed to confirm the theoretical results.

\section{Preliminaries and abstract framework}
\label{sec:Preliminaries}

Given two separable $\mathbb{R}$-Hilbert spaces $(U, \langle \cdot,\cdot \rangle_U,\|\cdot\|_U)$  and $(H, \langle \cdot,\cdot \rangle_H, \|\cdot\|_H)$, we denote by $\mathcal{L}(U, H)$ the Banach space of all linear bounded operators from $U$ into $H$ and $\mathcal{L}_2(U, H)$  the Hilbert space of all Hilbert-Schmidt operators from $U$ into $H$.
For simplicity, if $H = U$, we write $\mathcal{L}(U) =\mathcal{ L}(U,U)$
and $\mathrm{HS} = \mathcal{L}_2(U,U)$.  It is well known that
$
\|ST\|_{\mathcal{L}_2(U,H)}
\leq
 \|T\|_{\mathcal{L}_2(U,H)}
 \|S\|_{\mathcal{L}(U)},
\;
 \text{for}\;
 T\in \mathcal{L}_2(U,H),
 \;
 S\in \mathcal{L}(U).
$
Let $Q \in\mathcal{L}(U)$ be a self-adjoint, positive semidefinite operator. We denote the space of the Hilbert-Schmidt operators from
$Q^{\frac{1}{2}}(U)$ to $H$ by
$\mathcal{L}_2^0:
=
\mathcal{L}_2 (Q^{\frac{1}{2}}(U),H)$
and the corresponding
norm is given by
$
\|\Gamma \|_{\mathcal{L}_2^0}
=
\|\Gamma Q^{\frac{1}{2}}\|_{\mathcal{L}_2(U, H)}.$

Let $\textbf{E}$ be the expectation in the probability space and let $L^2(\Omega;H)$ be
the space of $H$-valued integrable random variables, equipped with the norm
$
\|v\|_{L^2(\Omega;H)}
:=
\big ( \mathbf{E} \big [ \|v\|^2 \big] \big) ^{\frac12}.
$
Next,  we introduce a self-adjoint, positive definite, linear operator $A = - L$ with the domain 
$D(A) = H^2 (\mathcal{D}) \cap H_0^1 (\mathcal{D})$. Then we define the separable 
Hilbert space $\dot{H}^s=D(A^{\frac{s}{2}})$, equipped with the corresponding norm
\begin{equation}
\|v\|_s
:=
\|A^{\frac s2} v\|
=
\Big(\sum_{j=1}^\infty\lambda_j^s \left <v,\phi_j\right>^2\Big)^{\frac12},
\;
s\in \mathbb{R},
\end{equation}
where $\{(\lambda_j,\phi_j)\}_{j=1}^\infty$ are the eigenpairs of $A$ with orthonormal eigenvectors.
It is well-known that $\dot{H}^0=H:=L_2(\mathcal{D})$, $\dot{H}^1=H_0^1(\mathcal{D})$ and $\dot{H}^2=H^2(\mathcal{D})\cap H_0^1(\mathcal{D})$.
To define the mild solution of \eqref{eq:stocha-str-dam-wave-equation} in the semigroup framework
as in \cite{da2014stochastic}, we additionally introduce the time derivative of  the solution $u$ as 
a new variable $v := u_t $ and rewrite (\ref{eq:stocha-str-dam-wave-equation}) as
\begin{align}\label{eq:st-dam-wave-abstract-form}
\left\{\begin{array}{ll}
\dd X(t)
=
- \mathcal{A}X(t)
\,
\dd t
+
\mathbf{F} (X(t)) \,
 \dd t
 +
  \mathbf{B}
  \,
  \dd W(t), &
t \in (0, T],\\
X(0)
=
X_0,&
\end{array}\right.
\end{align}
where $X(t):=(u(t),v(t))'$, $X_0:=(u_0,v_0)'$ and
\begin{equation}
\mathcal{A}
:=
\biggl[
\begin{array}{c
c}0&-I\\A&\alpha A\end{array}
\biggr],
\quad
\mathbf{F}(X)
:=
\biggl[
      \begin{array}{c}0 \\ F(u)\end{array}
\biggr]
\;\: \text{ and }
 \;
\mathbf{B}
:=
\biggl[
\begin{array}{c}0\\I\end{array}
\biggr].
\end{equation}
It was shown in \cite[Lemma 2.1]{larsson1991finite} that $-\mathcal{A}$ generates an analytic
semigroup $\mathcal{S}(t)$ in $\dot{H}^s\times \dot{H}^{s-\sigma}$,
for $ s \in \mathbb{R}, \sigma \in [0, 2]$, given by
\begin{align}\label{eq:definiton-semigroup-A}
\mathcal{S}(t)
=
e^{ - t\mathcal{A}}
=
\biggl[
\begin{array}{c c}
\mathcal{S}_1(t)
&
\mathcal{S}_2(t)
\\
\mathcal{S}_3(t)
&
\mathcal{S}_4(t)
 \end{array}
 \biggr].
\end{align}
Here, $\mathcal{S}_i(t)$, $i=1,2,3,4$,
can be expressed in terms of the eigenfunction expansion, for
$\varphi\in L_2(\mathcal{D})$,
\begin{align}\label{eq:definition-Si}
\left\{
\begin{array}{ll}
\mathcal{S}_1(t)\varphi
=
\sum_{j=1}^\infty
\frac{\lambda_j^+e^{-t\lambda_j^-}
         -\lambda_j^-e^{-t\lambda_j^+}}
             {\lambda_j^+-\lambda_j^-}
                 \left<\varphi,\phi_j\right>\phi_j,
                 &\\
\mathcal{S}_2(t)\varphi
=
\sum_{j=1}^\infty
              \frac{e^{-t\lambda_j^-}-e^{-t\lambda_j^+}}
                        {\lambda_j^+-\lambda_j^-}\left<\varphi,\phi_j\right>\phi_j,
                        &\\
\mathcal{S}_3(t)\varphi
=
\sum_{j=1}^\infty
             \frac{\lambda_j(e^{-t\lambda_j^+}-e^{-t\lambda_j^-})}
                         {\lambda_j^+-\lambda_j^-}\left<\varphi,\phi_j\right>\phi_j,
                         &\\
\mathcal{S}_4(t)\varphi
=
\sum_{j=1}^\infty
              \frac{\lambda_j^+e^{-t\lambda_j^+}-\lambda_j^-e^{-t\lambda_j^-}}
                      {\lambda_j^+-\lambda_j^-}\left<\varphi,\phi_j\right>\phi_j,
                      &
\end{array}\right.
\end{align}
  where $\{(\lambda_j,\phi_j)\}_{j=1}^\infty$ are the eigenpairs of $A$
   with orthonormal eigenvectors,
   and where $\lambda_j^{\pm}$ are the solutions of the following equation
\begin{equation}
z^2
-\alpha
\lambda_j z
+
\lambda_j
=0.
\end{equation}
As $\lambda_j \rightarrow \infty $,
one can easily judge that
$\lambda_j^-
   \thicksim 1/\alpha,
     \lambda_j^+\thicksim
      \alpha \lambda$.

To get the existence, uniqueness and regularity properties of the mild solution of
\eqref{eq:stocha-str-dam-wave-equation} and also for the purpose of the error analysis of the numerical scheme,
we give some assumptions on the nonlinear term, the noise process and the initial data as follows.
\begin{assumption}\label{assum:wiener-process-condition}
(Q-Wiener process). Let $W(t)$ be a (possibly cylindrical) $Q$-Wiener process,
with the covariance operator $Q \colon H \rightarrow H$ being a symmetric nonnegative operator satisfying
\begin{align}\label{jiezz1}
\|A^{\frac{\gamma-1}{2}}Q^{\frac{1}{2}}\|_{\mathrm{HS}}
<\infty \quad \text{ for some } \gamma\in(-1,2].
\end{align}
\end{assumption}
\begin{assumption} (Nonlinearity). \label{assum:F-condition}
The deterministic mapping $F \colon H \rightarrow H$ is assumed to be twice differentiable and  
there exists a positive constant $K$ such that, for $\beta := \min\{1,1+\gamma\}$,
\begin{align}
\|F(x)-F(y)\|
\leq
K \| x - y \|,
& \quad
\forall x,\;y
\in H,
\label{eq:F-assum-condition1}
\\
\|A^{-\frac{\theta}{2}}F'(y)z\|
\leq
 K (\|y\|_\beta+1) \|z\|_{-\beta},
 &
 \quad \forall y\in \dot{H}^\beta ,
 \;
  \forall z\in \dot{H}^{-\beta},
  \;
  \theta\in[1,2).
  \label{eq:F-assum-condition3}
\end{align}
\end{assumption}
In view of \eqref{eq:F-assum-condition1} in Assumption \ref{assum:F-condition}, the mapping $F$ also obeys
the following linear growth condition,
\begin{align}
\|F(x)\|
\leq
 K \|x\| + \| F (0) \|,
 \;
 \forall x\in H.
 \label{eq:growth-condtion-F}
\end{align}
Also, we remark that such condition as \eqref{eq:F-assum-condition3} was also used 
in \cite{wang2015note}, where the condition is validated only for particular ranges of 
$\beta$ as $\beta < \tfrac12$ and $\beta = 1$. In Example \ref{example:nonlinearity},
we give a class of concrete nonlinear Nemytskij operators to validate 
\eqref{eq:F-assum-condition3} for the whole range $\beta \in (0, 1]$.

\begin{assumption}(Initial data).\label{assum:initial-value}
 Let $u_0,v_0$ be $\mathcal{F}_0$-measurable and $(u_0,v_0)' \in 
 L^p(\Omega; \dot{H}^{\min\{2,\gamma+1\}}\times \dot{H}^\gamma)$, for any $p\geq 2$.
\end{assumption}

Under these assumptions, as shown in \cite[Theorem 2.1]{qi2015error},
the equation \eqref{eq:st-dam-wave-abstract-form}
has a unique  mild solution. Moreover, a slight modification of \cite[Theorem 2.1]{qi2015error} ensures
the following regularity results.
\begin{theorem}\label{thm:main-regularity}
Assume that Assumptions \ref{assum:wiener-process-condition}-\ref{assum:initial-value} hold. Then there exists a unique mild solution of \eqref{eq:st-dam-wave-abstract-form} given by
\begin{align}\label{eq:mild-SPDE}
X(t)=\mathcal{S}(t)X_0+\int_0^t\mathcal{S}(t-s)\mathbf{F}(X(s))\, \mathrm{d} s
+\int_0^t\mathcal{S}(t-s)\mathbf{B}\, \mathrm{d} W(s).
\end{align}
Furthermore, it holds that, for $p\geq 2$ and $t,s \in [0, T]$,
\begin{align}
\| u(t) \|_{ L^p( \Omega;\dot{H}^{\varrho} ) }+\| v(t) \|_{ L^p( \Omega;\dot{H}^{\varrho-2} ) }
       &\leq C \big( 1 + \|u_0\|_{L^p (\Omega; \dot{H}^{\varrho}) }
          + \|v_0\|_{L^p(\Omega; \dot{H}^{\varrho-2})}
           \big),
\label{eq:thm1-u-spatial2} \\
\| u(t) - u(s) \|_{L^p(\Omega;\dot{H}^0 )}
       \leq C |t-s|^{\frac{\varrho } 2 }
         & \big(1 + \|u_0\|_{L^p(\Omega; \dot{H}^{\varrho})}
          +\|v_0\|_{L^p(\Omega; \dot{H}^{\varrho-2})}
          \big),
 \label{eq:thm1-u-temporal2}
\end{align}
where we write  $ \varrho :=\min\{ 1+ \gamma, 2\}$. For the particular case when 
Assumption \ref{assum:wiener-process-condition} is fulfilled with $\gamma\in[0,2]$, 
it holds that, for $p\geq2$ and $t,s \in [0, T]$,
\begin{align}
\| u(t) \|_{L^p(\Omega;\dot{H}^{\gamma})}+\| v(t) \|_{L^p(\Omega;\dot{H}^{\gamma})}
      & \leq C \big( 1 + \|u_0\|_{L^p(\Omega; \dot{H}^{ \gamma})}
           + \|v_0\|_{L^p(\Omega; \dot{H}^{\gamma})}
           \big),
\label{eq:thm1-v-spatial2}\\
\| v(t) - v(s) \|_{L^p(\Omega;\dot{H}^{0})}
       \leq C |t-s|^{\frac{\min\{\gamma,1\}} {2} }
          &\big(1 + \|u_0\|_{L^p(\Omega; \dot{H}^{\min\{\gamma,1\}})}
          +\|v_0\|_{L^p(\Omega; \dot{H}^{\min\{\gamma,1\}})}
           \big).
 \label{eq:thm1-v-temporal2}
\end{align}
\end{theorem}
Here and below, $C$ represents a generic constant that may change between occurrences 
but only depends on $K, T, \| F (0) \|, \|A^{\frac{\gamma-1}{2}}Q^{\frac{1}{2}}\|_{\mathrm{HS}}$ and 
the initial data. As an immediate consequence, Theorem \ref{thm:main-regularity} implies the following facts.
\begin{lemma}
Under the assumptions of Theorem \ref{thm:main-regularity}, it holds that, with $ \varrho :=\min\{1+ \gamma, 2\}$,
\begin{align}
\sup_{ s \in [0, T] } \|F(u (s) )\|_{L^p(\Omega;\dot{H}^0)}
&
\leq
C\big(1 + \|u_0\|_{L^p(\Omega; \dot{H}^0)}
          +
          \|v_0\|_{L^p(\Omega; \dot{H}^{-2})}),
\label{lem:spatial-regularity-F}
\\
\|F(u(t))-F(u(s))\|_{L^p(\Omega;\dot{H}^0)}
&
\leq
      C |t-s|^{\frac{\varrho} {2} }
         \big(
         1 + \|u_0\|_{L^p(\Omega; \dot{H}^{\varrho})}
          +
          \|v_0\|_{L^p(\Omega; \dot{H}^{\varrho-2})}
          \big).
 \label{lem:temporal-regularity-F}
\end{align}
\end{lemma}
Subsequently, we present a class of concrete nonlinear mappings $F$ that fulfill Assumption \ref{assum:F-condition}.
\begin{example} \label{example:nonlinearity}
Let $H := L_2(\mathcal{D})$, $A=-\Delta=-\sum_{i=1}^d\frac{\partial^2}{\partial x_i^2}$ with the domain 
$D(A)=H^2(\mathcal{D})\cap H_0^1(\mathcal{D})$
and let $F:H\rightarrow H$ be a Nemytskij operator associated to $f$, defined by
\begin{equation}\label{ABSTAR1}
F(\varphi)(x)=f(x,\varphi(x)),\; x\in \mathcal{D},\; \varphi\in H,
\end{equation}
where $\mathcal{D}\subset \mathbb{R}^d$, $d=1,2,3$, is a bounded open set with Lipschitz boundary
and $f: \mathcal{D}\times \mathbb{R}\rightarrow \mathbb{R}$ is assumed to be a smooth nonlinear 
function satisfying, with some constant $K > 0$,
\begin{eqnarray}
&|f(x,z)|\leq K( |z| + 1 ),
&
\label{fassum1}\\
\Big|
\tfrac{\partial f}{\partial z}(x,z)
      \Big|
\leq
K,\;
&
\Big|
\tfrac{\partial^2 f}{\partial x_i\partial z}
(x,z)
     \Big|
     \leq
      K,
      &
       and\qquad
      \Big|
\tfrac{\partial^2 f} {\partial z^2}(x,z)
             \Big|
\leq K,
\label{fassum2}
\end{eqnarray}
for all $x=(x_1,x_2,\cdots, x_d)' \in \mathcal{D}$, $z\in \mathbb{R}$.
It is clear that $F$  satisfies the following conditions
\begin{eqnarray}
\|F(\varphi)\|
\leq
\sqrt{2} K ( \|\varphi\|+1 ),
\label{ABSTAR4}
\\
\|F(\varphi_1)-F(\varphi_2)\|
\leq
 K \|\varphi_1-\varphi_2\|,
    \label{ABSTAR5}
\end{eqnarray}
for $\varphi, \varphi_1, \varphi_2\in H$. The first derivative operator of $F$ is given by
\begin{eqnarray}
F'(\varphi)(\psi)(x)
=
\tfrac{\partial f}{\partial z}
           (x,\varphi(x))
               \cdot \psi(x),\;
 x\in \mathcal{D},\label{ABSTAR2}
\end{eqnarray}
for all $\varphi,\psi\in H$. Note that the derivative operator $F'(\varphi),  \varphi\in H$
defined in the above way is self-adjoint and
 by similar arguments used in the proof of \cite[Lemma 1]{wang2014higher},
 we have, for $\beta=\min\{1,1+\gamma\}\in (0,1]$,
\begin{equation}\label{spr4}
\|F'(\varphi)\psi\|_{W^{\beta,2}}
\leq
 C(\|\psi(x)\|+\|\psi(x)\|_{C(\mathcal{D},\mathbb{R})}
 +
 \|\varphi(x)\|_{W^{\beta,2}}\|\psi(x)\|_{C(\mathcal{D},\mathbb{R})})
+
C\|\psi(x)\|_{W^{\beta,2}},
\end{equation}
where $\varphi\in \dot{H}^\beta$, $\psi\in \dot{H}^\theta$ with $\theta\in [1,2)$ and 
$\theta >\frac{d}{2}$, $d=1,2,3$.  Here $ \| \cdot \|_{W^{\beta,2}} $ stands for the usual
norm for the Sobolev space $W^{\beta,2}$. It is well-known that
\begin{eqnarray}\label{eq:space-norm-relation}
\frac{1}{C_\beta}\|v\|_{W^{\beta,2}}
\leq
\|v\|_\beta
\leq
C_\beta\|v\|_{W^{\beta,2}},
 \;
 \text{ for } \; v\in \dot{H}^\beta.
\end{eqnarray}
Furthermore, the Sobolev embedding theorem gives $\dot{H}^\theta\subset C(\mathcal{D},\mathbb{R})$ 
and $\dot{H}^\theta \subset \dot{H}^\beta$, for $\theta>\frac d 2$ and  $\theta\in [1,2)$.
Due to \eqref{eq:space-norm-relation}, we infer that
\begin{eqnarray}\label{spr6}
\|F'(\varphi)\psi\|_{W^{\beta,2}}
\leq
C(\|\varphi\|_\beta+1)
                            \|\psi\|_\theta,
\;for\;
       \varphi\in \dot{H}^\beta,
       \;
       \psi\in \dot{H}^\theta.
\end{eqnarray}
 For $\varphi\in \dot{H}^\beta$ and $\psi\in \dot{H}^\theta\subset H_0^1(\mathcal{D})$, we have $F'(\varphi)\psi\in \dot{H}^\beta$ and
\begin{eqnarray}
\|F'(\varphi)\psi\|_\beta
\leq
\frac 1 C_\beta
           \|F'(\varphi)\psi\|_{W^{\beta,2}}
\leq
C(\|\varphi\|_\beta+1)
                  \|\psi\|_\theta.
\end{eqnarray}
Consequently, one can show that, for $\varphi\in \dot{H}^\beta$ and $\psi\in \dot{H}^{-\beta}$,
\begin{align}
\|A^{-\frac{\theta}{2}}F'(\varphi)\psi\|
&=
\sup_{\chi\in \dot{H}^0}
            \frac{\big< A^{-\frac{\theta}2}F'(\varphi)\psi,\chi\big> }{ \|\chi\| }
=
\sup_{\chi\in \dot{H}^0}
              \frac{\big<A^{-\frac\beta2}\psi,A^{\frac\beta2}F'(\varphi)A^{-\frac{\theta}{2}}\chi\big>}
              {\|\chi\|}
\nonumber\\
&\leq
\sup_{\chi\in \dot{H}^0}
           \frac{\|A^{-\frac\beta2}\psi\|\|A^{\frac\beta2}F'(\varphi)A^{-\frac\theta2}\chi\|}{\|\chi\|}
\leq
 \sup_{\chi\in \dot{H}^0}
               \frac{\|\psi\|_{-\beta}
                           \|A^{\frac\beta2}F'(\varphi)A^{-\frac\theta2}\chi\|}{\|\chi\|}
 \nonumber\\&\leq
 C \sup_{\chi\in \dot{H}^0}
              \frac{\|\psi\|_{-\beta}\left(1+\|\varphi\|_\beta\right)\|\chi\|}{\|\chi\|}
  \nonumber\\
  &\leq
  C\|\psi\|_{-\beta}(1+\|\varphi\|_\beta).
\end{align}
To conclude, the function $F$ defined by \eqref{ABSTAR1} satisfies Assumption \ref{assum:F-condition}.
\end{example}

\section{The proposed scheme and main result}
\label{sec:main-result}

In this section, we construct a full discrete scheme for \eqref{eq:stocha-str-dam-wave-equation},
with the spatial approximation done by a spectral Galerkin method,
along with temporal discretization by an accelerated exponential Euler scheme. To do this, we define a finite dimensional subspace of $H$ by $H_N:= \text{span} \{\phi_1,\phi_2,\cdots,\phi_N\}$
and the projection operator $\mathcal{P}_N: \dot{H}^\kappa \rightarrow H_N$ by
\begin{align}
\mathcal{P}_Nv=\sum_{j=1}^N\left<v,\phi_j\right>\phi_j,\;\forall v\in \dot{H}^\kappa,\; \kappa \in \mathbb{R}.
\end{align}
Then one can straightforwardly show that $\|\mathcal{P}_Nv\|\leq \|v\|$ and that, for $v\in \dot{H}^\kappa,\;\forall\kappa\geq 0$
\begin{align}\label{eq:interpolation-property-PN}
\|(I-\mathcal{P}_N)v\|
=
\Big(
\sum_{j={N+1}}^\infty\left<v,\phi_j\right>^2
      \Big)^{\frac12}
\leq
\lambda_{N+1}^{-\frac\kappa2}
\Big(
\sum_{j={N+1}}^\infty\lambda_j^\kappa\left<v,\phi_j\right>^2
       \Big)^{\frac12}
\leq
 \lambda_{N+1}^{-\frac\kappa2}
    \|v\|_\kappa.
\end{align}

Let $A_N \colon H_N \rightarrow H_N$ be defined as $A_N=\mathcal{P}_NA$. Then a spatial approximation of
\eqref{eq:st-dam-wave-abstract-form} leads to a stochastic differential equation in $H_N\times H_N$,
\begin{align}
\label{eq:st-spactial-dam-wave-abstract-form}
\left\{
\begin{array}{ll}
\dd X^N(t)
=
- \mathcal{A}_NX^N(t)
 \, \dd t
+
\mathbf{F}_N (X^N(t))
 \, \dd t
 +
  \mathbf{B}_N
  \, \dd W(t),
t \in (0, T],
 \\
   X^N(0)
   =
   X^N_0,
\end{array}
\right.
\end{align}
where $X^N(t):=(u^N(t),v^N(t))'$,
   $X^N_0:=(\mathcal{P}_Nu_0,\mathcal{P}_Nv_0)'$
     and
\[
\mathcal{A}_N
:=
\biggl[
\begin{array}{c c}
       0&-I\\A_N&\alpha A_N
\end{array}
      \biggr],
      \quad
\mathbf{F}_N(X^N)
:=
\biggl[
      \begin{array}{c}0 \\ P_NF(u^N)\end{array}
       \biggr]
       \;\:
       \text{ and }
\;
\mathbf{B}_N
:=
\biggl[
    \begin{array}{c}0\\P_N
       \end{array}
       \biggr].
\]
Similarly, $\mathcal{A}_N$ also generates an analytic semigroup $\mathcal{S}_N(t)$, $t\geq 0$,
in $H_N\times H_N$ which is an analogue of \eqref{eq:definiton-semigroup-A}, expressed by
 \begin{align}\label{eq:definiton-semigroup-AN}
\mathcal{S}_N(t)
=
e^{t\mathcal{A}_N}
=
\biggl[
\begin{array}{c c}
    \mathcal{S}_{1,N}(t)&\mathcal{S}_{2,N}(t)
    \\
    \mathcal{S}_{3,N}(t)&\mathcal{S}_{4,N}(t)
    \end{array}
    \biggr],
\end{align}
where $\mathcal{S}_{i,N}(t)$, $i=1,2,3,4$ are defined as follows,  for $\varphi \in H_N$,
\begin{align}\label{eq:definition-SiN}
\left\{
   \begin{array}{ll}
\mathcal{S}_{1,N}(t)\varphi
=
\sum_{j=1}^N
     \frac{\lambda_j^+e^{-t\lambda_j^-}-\lambda_j^-e^{-t\lambda_j^+}}
                                                      {\lambda_j^+-\lambda_j^-}
                                                      \left<\varphi,\phi_j\right>\phi_j,
     &\\
\mathcal{S}_{2,N}(t)\varphi
=
\sum_{j=1}^N
      \frac{e^{-t\lambda_j^-}-e^{-t\lambda_j^+}}
                                {\lambda_j^+-\lambda_j^-}
                                \left<\varphi,\phi_j\right>\phi_j,
&\\
\mathcal{S}_{3,N}(t)\varphi
=
 \sum_{j=1}^N
         \frac{\lambda_j(e^{-t\lambda_j^+}-e^{-t\lambda_j^-})}
                                      {\lambda_j^+-\lambda_j^-}
                                      \left<\varphi,\phi_j\right>\phi_j,
 &\\
\mathcal{S}_{4,N}(t)\varphi
=
\sum_{j=1}^N
      \frac{\lambda_j^+e^{-t\lambda_j^+}-\lambda_j^-e^{-t\lambda_j^-}}
                                      {\lambda_j^+-\lambda_j^-}
                                      \left<\varphi,\phi_j\right>\phi_j.&
\end{array}\right.
\end{align}
Since $A$ commutes with $P_N$,  Assumptions \ref{assum:wiener-process-condition}-\ref{assum:initial-value}
suffice to guaranteen  a unique solution $X^N(t)=(u^N(t),v^N(t))'$ for \eqref{eq:st-spactial-dam-wave-abstract-form}, given by
\begin{align}\label{eq:mild-spatial-SPDE}
X^N(t)
=
\mathcal{S}_N(t)X^N_0
+
\int_0^t\mathcal{S}_N(t-s)\mathbf{F}_N(X^N(s))\, \dd s
+
\int_0^t\mathcal{S}_N(t-s)\mathbf{B}_N\, \dd W(s).
\end{align}

Next we turn our attention to the full discretization of (\ref{eq:st-dam-wave-abstract-form}).  Given a time step-size
$k=T/M$ for some $M\in \mathbb{N}^+$,  we construct uniform mesh grids $t_m=mk$, for $m=0,1,2,\cdots M$. On the grids, we introduce the following full-discrete scheme based on the spatial discretization \eqref{eq:mild-spatial-SPDE}:
\begin{align}\label{eq:full-discretization-scheme-abstract}
X^N_{m+1}
=
\mathcal{S}_N(k)X^N_{m}
+
k\mathcal{S}_N(k)\mathbf{F}_N(X_m^N)
+
\int_{t_{m}}^{t_{m+1}}\mathcal{S}_N(t_{m+1}-s)\mathbf{B}_N\,\dd W ( s ).
\end{align}
With the notation in \eqref{eq:definition-SiN}, one can rewrite (\ref{eq:full-discretization-scheme-abstract}) as
\begin{align}\label{eq:full-scheme-form}
\left\{\begin{array}{ll}
u_{m+1}^N
=
\mathcal{S}_{1,N}(k)u_m^N
+
\mathcal{S}_{2,N}(k)v_m^N
+
k\mathcal{S}_{2,N}(k)P_NF(u_m^N)
+
\int_{t_m}^{t_{m+1}}\mathcal{S}_{2,N}(t_{m+1}-s)\,\dd W(s),
&\\
v_{m+1}^N
=
\mathcal{S}_{3,N}(k)u_m^N
+
\mathcal{S}_{4,N}(k)v_m^N
+
k\mathcal{S}_{4,N}(k)P_NF(u_m^N)
+
\int_{t_m}^{t_{m+1}}\mathcal{S}_{4,N}(t_{m+1}-s)\,\dd W(s),
&
\end{array}\right.
\end{align}
for $m=0,1,2,\cdots,M-1$, $M\in \mathbb{N}^+$.
%
%
%
A natural question arises how to simulate this scheme. More precisely, the key point is how to simulate
the two stochastic convolutions  $\int_{t_m}^{t_{m+1}}\mathcal{S}_{2,N}(t_{m+1}-s)\,\dd W(s)$
and $\int_{t_m}^{t_{m+1}}\mathcal{S}_{4,N}(t_{m+1}-s)\,\dd W(s)$ in \eqref{eq:full-scheme-form}.
It is not difficult to find that they are two correlated normally distributed random variables.
%
Once the operators $A$ and $Q$ own the same eigenfunctions, the correlation between the stochastic
convolutions can be explicitly computed and  the proposed scheme becomes rather easy to implement.
To show this, from here and below in this section we assume
\begin{equation}
Q \phi_j = \gamma_j \phi_j, \quad j \in \N,
\end{equation}
where $\{(\lambda_j,\phi_j)\}_{j=1}^\infty$ are the eigenpairs of $A$.
For $m=0,1,2,\cdots,M-1$, $j=1,2,\cdots, N$,
the following  random variables
\begin{align}
\eta_m^j
:=
\Big<\int_{t_m}^{t_{m+1}}
      \mathcal{S}_{2,N}(t_{m+1}-s)\,\dd W(s),\phi_j
      \Big>
=
\int_{t_m}^{t_{m+1}}
\frac{e^{-(t_{m+1}-s)\lambda_j^-}-e^{-(t_{m+1}-s)\lambda_j^+}}{\lambda_j^+ - \lambda_j^- }
\gamma_j^{\frac{1}{2}} \,\dd \beta_j(s),\nonumber
\end{align}
and
\begin{align}
\widehat{\eta}_m^j
:=
\Big<
\int_{t_m}^{t_{m+1}}\mathcal{S}_{4,N}(t_{m+1}-s)\dd W(s),\phi_j
       \Big>
=
\int_{t_m}^{t_{m+1}}
                  \frac{ \lambda_j^+e^{-(t_{m+1}-s)\lambda_j^+ }-\lambda_j^-e^{-(t_{m+1}-s)\lambda_j^-}  }
                     {\lambda_j^+-\lambda_j^-}
\gamma_j^{\frac12} \dd \beta_j(s),
\nonumber
\end{align}
are two series of mutually independent, normally distributed random variables, 
where $\lambda_j^+$ and $\lambda_j^-$ are defined as in section \ref{sec:Preliminaries}.
Furthermore, they satisfy $\mathbf{E}[\eta_m^j]=0$,  $\mathbf{E}[\widehat{\eta}_m^j]=0$,
\begin{align}
\mathrm{Var}(\eta_m^j)
=
\mathbf{E}[|\eta_m^j|^2]
=
\frac{\gamma_j}{(\lambda_j^+-\lambda_j^-)^2}
\Big[
       \frac \alpha 2 - \frac {e^{-2k\lambda_j^+}}  {2\lambda_j^+}
                 -\frac {e^{-2k\lambda_j^-}}  {2\lambda_j^-}
                      -\frac{2}{\alpha}(1-e^{-2k\lambda_j})
\Big],
\end{align}
and
\begin{align}
\mathrm{Var}(\widehat{\eta}_m^j)
=
\mathbf{E}[|\widehat{\eta}_m^j|^2]
=
\frac{\gamma_j}{(\lambda_j^+-\lambda_j^-)^2}
&\Big[
\frac{\lambda_j}{2}
                      -\frac{\lambda_j^+}{2}e^{-2k\lambda_j^+}
           -\frac{\lambda_j^-}{2}e^{-2k\lambda_j^-}
                                                  -\frac{2}{\alpha\lambda_j}(1-e^{-\alpha\lambda_j k})
\Big].
\end{align}
Let $M_m^j$ be a family of $2\times 2$ matrices with
\begin{align}
C_m^i:=\biggl[
           \begin{array}{c c}
            \mathrm{Var}(\eta_m^j)
&
\mathrm{Cov}(\eta_m^j,\widehat{\eta}_m^j)
\\
\mathrm{Cov}(\eta_m^j,\widehat{\eta}_m^j)
&
\mathrm{Var}(\widehat{\eta}_m^j)
\end{array}
\biggr]
=
M_m^j(M_m^j)^T,
\end{align}
where $\mathrm{Cov}(\eta_m^j,\widehat{\eta}_m^j)$ are the covariances of $\eta_m^j$ and $\widehat{\eta}_m^j$ given by
\begin{align}
\mathrm{Cov}(\eta_m^j,\widehat{\eta}_m^j)
=
\mathbf{E}[\eta_m^j\widehat{\eta}_m^j]
=
\frac {\gamma_j} {(\lambda_j^+-\lambda_j^-)^2}
                                  \Big(\frac {e^{-2k\lambda_j^+}+e^{-2k\lambda_j^-}} 2
           -e^{-\alpha k\lambda_j}
                                \Big).
\end{align}
Hence the pair of random variables $(\eta_m^j, \widehat{\eta}_m^j)$ can be fully characterized by
\begin{align}
&\biggl[
\begin{array}{c}\eta_m^j
\\
\widehat{\eta}_m^j
\end{array}
\biggr]
=
M_m^j
\biggl[
\begin{array}{c}
\xi_m^j\\
\widehat{\xi}_m^j
\end{array}
\biggr],
\end{align}
with $\xi_m^j$ and $\widehat{\xi}_m^j$, for $j=1,2,\cdots, N$ and $m=0,1,2,\cdots, M-1$ being 
independent, standard normally distributed random variables.
Then for $j=1,2,\cdots,N$ and $m=0,1,2,\cdots, M-1$ the coefficients of the expansion of $u_m^N$ and $v_m^N$ in the scheme \eqref{eq:full-scheme-form} can be realized by the following recurrence equation:
\begin{align}
\left<u_{m+1}^N,\phi_j\right>
=
\frac{\lambda_j^+e^{-t\lambda_j^-}
-\lambda_j^-e^{-k\lambda_j^+}}{\lambda_j^+-\lambda_j^-}\left<u_m^N,\phi_j\right>
&+
\frac{e^{-k\lambda_j^-}-e^{-k\lambda^+_j}}
                                                       {\lambda_j^+-\lambda_j^-}
\left<v_m^N,\phi_j\right>
\nonumber\\
&
+
k\frac{e^{-k\lambda_j^-}-e^{-k\lambda_j^+}}{\lambda_j^+-\lambda_j^-}
\left<F(u_m^N),\phi_j\right>+\eta_m^j,
\\
\left<v_{m+1}^N,\phi_j\right>
=
\frac{\lambda_j(e^{-t\lambda_j^+}
                               -e^{-k\lambda_j^-})}{\lambda_j^+-\lambda_j^-}
\left<u_m^N,\phi_j\right>
&+
\frac{\lambda_j^+e^{-k\lambda_j^+}-\lambda_j^-e^{-k\lambda_j^-}}
                                                 {\lambda_j^+-\lambda_j^-}
\left<v_m^N,\phi_j\right>
\nonumber\\
&+
k\frac{\lambda_j^+e^{-k\lambda_j^+}-\lambda_j^-e^{-k\lambda_j^-}}
                                   {\lambda_j^+-\lambda_j^-}
\left<F(u_m^N),\phi_j\right>+\widehat{\eta}_m^j.
\end{align}

Now we state our main results.
\begin{theorem}\label{them:main-results}
Suppose that Assumptions \ref{assum:wiener-process-condition}-\ref{assum:initial-value} hold.
Let $(u(t),v(t))'$ be the solution of \eqref{eq:st-dam-wave-abstract-form} given by \eqref{eq:mild-SPDE}
and $(u_m^N,v_m^N)'$ be the numerical solution produced by  \eqref{eq:full-scheme-form}. 
Then it holds that, for $\gamma\in (-1,2]$,
\begin{align}\label{them:error-estimates-u-um}
\|u(t_m)-u^N_m\|_{L^2(\Omega;\dot{H}^0)}
\leq
  C(k^{\min\{1+\gamma,1\}}+\lambda_{N+1}^{-\frac{1+\min\{\gamma, 1\}}{2}}),
\end{align}
and for  $\gamma\in [0,2]$,
\begin{align}\label{them:error-estimates-v-vm}
\|v(t_m)-v^N_m\|_{L^2(\Omega;\dot{H}^0)}
\leq
 C(k+\lambda_{N+1}^{-\frac{\gamma}{2}}).
\end{align}
\end{theorem}

Theorem \ref{them:main-results} reveals that the obtained strong orders in space are optimal in the sense that the convergence orders coincide with the orders of the spatial regularity of $(u,v)'$
 as stated in Theorem \ref{thm:main-regularity}. However, the strong order
in time is twice as high as the order of the time regularity of the mild solution for 
$\gamma\in (-1,0]$, especially for the case of the space-time white noise in two or three dimensions.
For  $\gamma\in[0,2]$, $v_m^N$ as well as  $u_m^N$  strongly converges with a rate of order one, 
regardless of the order of regularity of the velocity $v$.  As already discussed in the introduction,  
this scheme allows for higher strong order in time than the linear implicit
 Euler scheme investigated in \cite{qi2015error}.

\section{Proof of the main result}
\label{sec:proof}
This section is devoted to the proof of the main convergence result Theorem  \ref{them:main-results}.
In the first part, we present some useful preparatory results.
In particular, the following proposition plays an essential role in the error analysis of our scheme.
\begin{proposition}\label{pro:temporal-regularity}
Under the assumptions of Theorem \ref{thm:main-regularity}, it holds that, for $\gamma\in (-1,2]$,
\begin{align}
\| u(t) - u(s) \|_{L^p(\Omega;\dot{H}^{-\delta})}
       \leq C |t-s|^{\frac{\varrho +\delta} {2} }
         & \big(1 + \|u_0\|_{L^p(\Omega; \dot{H}^{\varrho})}
          +\|v_0\|_{L^p(\Omega; \dot{H}^{\varrho-2})}
          \big),\quad t \geq s,
 \label{eq:thm1-u-temporal2}
\end{align}
where
\begin{align}\label{constant:order-u}
\delta \in [0, 1 - \min\{ 1, | \gamma | \} ],
\quad \mathrm{and} \quad
\varrho := \min\{\gamma + 1, 2\}.
%
\end{align}
\end{proposition}
Before showing Proposition \ref{pro:temporal-regularity},  we quote some useful results from \cite[Lemma 2.9]{qi2015error}. To do this,  we introduce two operators $P_1$ and $P_2$ defined by
\[
P_ix=x_i,\;\forall x=\left(x_1,x_2\right)'\in H_1\times H_2
\]
for two Hilbert spaces $H_1$, $H_2$.
Then $P_1 \mathcal{S}(t)x=\mathcal{S}_1(t)x_1+\mathcal{S}_2(t)x_2$ and $P_2 \mathcal{S}(t)x=\mathcal{S}_3(t)x_1+\mathcal{S}_4(t)x_2$.
\begin{lemma} \label{lem:existing-regularity}
Denote $x=(x_1,x_2)'$ and let  $\mathcal{S}(t)$ be the semigroup defined as above, with four components 
$\mathcal{S}_i(t)$, $i=1,2,3,4$  as in \eqref{eq:definition-Si}. Then for $0<s<t\leq T$ the following regularity 
properties hold true,
\begin{align}
\|P_1\mathcal{S}(t)x\|
&\leq
C(\|x_1\|
+
\|x_2\|_{-\mu}),
&
 \;
 \mu\in[0,2],
\label{lem:spatial-regularity-S1-S2}
\\
\|P_1[\mathcal{S}(t)-\mathcal{S}(s)]x\|
&\leq
C(t-s)^{\frac \mu 2}(\|x_1\|_\mu
+
\|x_2\|_{ \mu-2 } ),
&
\mu\in[0,2],
\label{lemma:temporal-regularity-S1-S2}
\\
\| P_2 [ \mathcal{S}(t)-\mathcal{S}(s)]x  \| & \leq
C(t-s)^{\frac \mu 2}(\|x_1\|_\mu
+
\|x_2\|_\mu  ),
&
\mu \in [0, 2],
\label{lemma:temporal-regularity-S3-S4}
\\
\int_s^t
\|P_2\mathcal{S}(t-r)x\|^2\,\mathrm{d} r
&\leq
 C(t-s)^{ \nu }
 (
\|x_1\|^2_{ \nu  }
 +
 \|x_2\|^2_{ \nu - 1 }
 ),
 &
  \;
 \nu\in[0,1].
 \label{lem:integral-regularity1-S4}
\end{align}
\end{lemma}
{ \it  Proof of Proposition \ref{pro:temporal-regularity}. }
Noting that $ \mathcal{S}(0) = I$, we employ \eqref{eq:definiton-semigroup-A} and \eqref{eq:mild-SPDE} to obtain
\begin{align}
u(t)-u(s)
=
P_1[\mathcal{S}(t-s)-\mathcal{S}(0)]X(s)
+
\int_s^t\mathcal{S}_2(t-r)F(u(r))\,\dd r
+
\int_s^t\mathcal{S}_2(t-r)\,\dd W(r),
\label{eq:difference-u}
\end{align}
which suggests that
\begin{align}\label{eq:error-u(t)-u(s)}
\|u(t)-u(s)\|_{L^p(\Omega;\dot{H}^{-\delta})}
\leq&
 \|P_1[\mathcal{S}(t-s)-\mathcal{S}(0)]X(s)\|_{L^p(\Omega; \dot{H}^{-\delta } )}
\nonumber\\&
+
\int_s^t\|\mathcal{S}_2(t-r)F(u(r))\|_{L^p(\Omega; \dot{H}^{-\delta })}\,\dd r
+
\Big\|\int_s^t\mathcal{S}_2(t-r)\,\dd W(r)\Big\|_{L^p(\Omega; \dot{H}^{-\delta } )}
\nonumber\\
:= & 
\mathbb{I}+\mathbb{II}+\mathbb{III}.
\end{align}
Since $\frac{\varrho+\delta}{2}\leq 1$, for $\delta,\; \varrho$ as defined in \eqref{constant:order-u},  one can derive  by \eqref{eq:thm1-u-spatial2} and \eqref{lemma:temporal-regularity-S1-S2},
\begin{align}\label{eqn:partI-error-u-uh}
\mathbb{I}
&=
\| P_1[\mathcal{S}(t-s)-\mathcal{S}(0)] \big(A^{-\frac{\delta}{2}} u (s), A^{-\frac{\delta}{2}} v (s) \big)' \|_{L^p(\Omega; \dot{H}^{0 } )}
\nonumber
\\
&\leq
C(t-s)^{\frac{\varrho+\delta}{2}}
(\|u(s)\|_{L^p(\Omega;\dot{H}^{\varrho })}
+
\|v(s)\|_{L^p(\Omega;\dot{H}^{\varrho-2})})\nonumber\\
&\leq
C(t-s)^{\frac{\varrho+\delta}{2}}(1+\|u_0\|_{L^p(\Omega;\dot{H}^{\varrho })}
+
\|v_0\|_{L^p(\Omega;\dot{H}^{\varrho-2})}).
\end{align}
For  $\mathbb{II}$, using \eqref{lem:spatial-regularity-F} and applying \eqref{lem:spatial-regularity-S1-S2} with $x=(0, F(u))'$ and $\mu=0$
enable us to obtain that
\begin{align}\label{eqn:partII-error-u-uh}
 \mathbb{II}
\leq
\int_s^t \|F(u(r))\|_{L^p(\Omega; \dot{H}^0)}\,\dd r
 \leq
C(t-s)(1+\|u_0\|_{L^p(\Omega;\dot{H}^{0})}
+
\|v_0\|_{L^p(\Omega;\dot{H}^{-2})}).
\end{align}
To handle the estimate of $\mathbb{III}$, we utilize \eqref{lemma:temporal-regularity-S1-S2}, the Burkholder-Davis-Gundy type inequality and the definition of the Hilbert-Schmidt norm  to get
\begin{align}\label{eqn:partIII-error-u-uh}
\mathbb{III}^2
&\leq
 C_p^2  \int_s^t\big\|\mathcal{S}_2(t-r)A^{-\frac \delta 2}Q^{\frac12}\big\|_{\mathrm{HS}}^2\,\dd r
 \nonumber\\
&=
C_p^2\sum_{j=1}^\infty\int_s^t
\big\|
[\mathcal{S}_2(t-r)-\mathcal{S}_2(0)]A^{-\frac \delta 2}Q^{\frac12}
\phi_j\big\|^2
\,\dd r
\nonumber\\
&
\leq
 C\sum_{j=1}^\infty
 \int_s^t
(t-r)^{\varrho+\delta}
\big\|
A^{\frac {\varrho-2} 2}Q^{\frac{1}{2}}
\phi_j
\big\|^2
\,\dd r\nonumber\\
&\leq
 C
(t-s)^{1+\varrho+\delta}
\big\|
A^{\frac {\varrho-2} 2}Q^{\frac12}
\big\|_{\mathrm{HS}}^2
\nonumber\\
&
\leq  C
(t-s)^{1+\varrho+\delta}
\big\|A^{\frac{\gamma-1}{2}}Q^{\frac{1}{2}}\big\|_{\mathrm{HS}}^2,
\end{align}
where we also used the fact $\mathcal{S}_2(0)x=0$. At last, inserting
\eqref{eqn:partI-error-u-uh}-\eqref{eqn:partIII-error-u-uh} into
 \eqref{eq:error-u(t)-u(s)}
finishes the proof of this proposition.  $\square$

To complete the proof of the main result, we additionally need to exploit further regularity results on the
semigroup $\mathcal{S}(t)$, as stated in the following lemmas.
\begin{lemma}\label{lem:spatial-regualrity-semigroup}
Assume that  $\mathcal{S}(t)$ is the analytic semigroup with four components $\mathcal{S}_i(t)$, $i=1,2,3,4$
as defined in \eqref{eq:definition-Si}. Then for $x=(x_1,x_2)'$ and $0<s<t\leq T$ we have
\begin{align}
\|P_2\mathcal{S}(t)x\|
&\leq
C(\|x_1\|
+
t^{-\frac \mu 2}\|x_2\|_{-\mu}),
&
\;
\mu\in[0,2],
\label{lem:spatial-regularity-S3-S4}
\\
\|P_2[\mathcal{S}(t)-\mathcal{S}(s)]x\|
&\leq
 C(t-s)^{\frac\mu2}(\|x_1\|_\mu
 +
 s^{-\frac{\mu-\nu}2}\|x_2\|_\nu),
 &
\;
\mu\in[0,2],\;\nu\in[-2,\mu],
\label{lemma:temporal-regularity-S4}
\\
\Big\|
\int_s^t
\mathcal{S}_4(t-r)x_2\,
\mathrm{d} r
       \Big\|
&\leq
     C(t-s)^{\frac{2-\mu}{2}}
           \|x_2\|_{-\mu},
&
     \; \mu\in[0,2].
   \label{lem:integral-regularity2-S4}
\end{align}
\end{lemma}
{ \it  Proof of Lemma \ref{lem:spatial-regualrity-semigroup}. }
Thanks to the interpolation theory,  for \eqref{lem:spatial-regularity-S3-S4} we only need to verify 
the two cases $\mu=0$ and $\mu=2$, which can be directly obtained 
from \cite[Lemma 2.3]{larsson1991finite}.

Concerning \eqref{lemma:temporal-regularity-S4},  we first use \eqref{lemma:temporal-regularity-S3-S4}
with the choice of $x = (x_1, 0)'$ results in the following estimates
\begin{align}
\|[\mathcal{S}_3(t)-\mathcal{S}_3(s)]x_1\|
\leq
 C(t-s)^{\frac{\mu}{2}}\|x_1\|_\mu,\;\forall x_1\in \dot{H}^{\mu},\;\mu\in [0,\; 2].
\end{align}
Therefore, bearing the definitions of $P_2$ and $\mathcal{S}(t)$ in mind,  it suffices to prove
\begin{align}
\|[\mathcal{S}_4(t)-\mathcal{S}_4(s)]x_2\|
\leq
 C (t-s)^{\frac{\mu}{2}} s^{-\frac{\mu-\nu}{2}}  \|x_2\|_\nu,
 \; \forall x_2\in \dot{H}^{\nu},
 \; \mu\in[0,2],\;\nu\in[-2,\mu].
\end{align}
Employing \eqref{eq:definition-Si} implies
\begin{align}\label{eq:decompose-S4}
\|[\mathcal{S}_4(t)-\mathcal{S}_4(s)]x_2\|^2
=&
\sum_{j=0}^\infty
    \frac{[\lambda^+_j(e^{-t\lambda^+_j}
            -e^{-s\lambda^+_j})-\lambda^-_j(e^{-t\lambda^-_j}
                   -e^{-s\lambda^-_j})]^2}{[\lambda^+_j-\lambda^-_j]^2}
\left<x_2,\phi_j\right>^2
\nonumber\\
\leq&
2\sum_{j=0}^\infty
               \frac{(\lambda^+_j)^2(e^{-t\lambda^+_j}
                    -e^{-s\lambda^+_j})^2}
                           {[\lambda^+_j-\lambda^-_j]^2}
\left
<x_2,\phi_j\right>^2
\nonumber\\
&+
  2\sum_{j=0}^\infty\frac{(\lambda^-_j)^2(e^{-t\lambda^-_j}
                     -e^{-s\lambda^-_j})^2}
                        {[\lambda^+_j-\lambda^-_j]^2}
\left<x_2,\phi_j\right>^2
\nonumber \\ := &
B_1 + B_2.
\end{align}
Before proceeding further,  we recall that
\begin{align}\label{eq:eqvlent}
\lambda_j^- \thicksim 1/\alpha, \;  \lambda_j^+ \thicksim \lambda_j, \; \text{ as } \; \lambda_j \rightarrow \infty.
\end{align}
With this and the fact that $\sup_{s>0}s^{\mu-\nu} e^{-s}<\infty$, for $\nu\leq \mu$, we start the estimate of $B_1$ as follows:
\begin{align}\label{eq:A-estimation}
B_1
&\leq
C\sum_{j=0}^\infty
      (e^{-t\lambda^+_j}-e^{-s\lambda^+_j})^2
                   \left<x_2,\phi_j\right>^2
\nonumber\\
&\leq
 C\sum_{j=0}^\infty
  \left[
            (\lambda^+_j)^{\frac{\mu-\nu}{2}}
                      e^{-s\lambda^+_j}
                           (\lambda^+_j)^{-\frac{\mu}{2}}
                            (e^{-(t-s)\lambda^+_j}-1)
                                                    \right]^2
(\lambda^+_j)^\nu
      \left<x_2,\phi_j\right>^2
\nonumber\\
&
\leq
    Cs^{-(\mu-\nu)}(t-s)^{\mu}
                    \|x_2\|_\nu^2,
\end{align}
where we also used $(\lambda^+_j)^{-\frac \mu 2}(1-e^{-(t-s)\lambda^+_j})\leq C(t-s)^{\frac\mu2}$, for $0\leq\mu\leq 2$.
In the same manner, we get
\[B_2
\leq
Cs^{-(\mu-\nu)}(t-s)^{\mu}
\|x_2\|_\nu^2.
\]
Inserting this estimate and \eqref{eq:A-estimation} into \eqref{eq:decompose-S4}  ends the  proof of \eqref{lemma:temporal-regularity-S4}.
With regard to \eqref{lem:integral-regularity2-S4}, the definition of $\mathcal{S}_4(t)$ enables us to deduce that
\begin{align}
\Big\|
\int_s^t\mathcal{S}_4(t-r)x_2\,\dd r
      \Big\|
=
\bigg(
   \sum_{j=1}^\infty
    \bigg|
    \int_s^t
        \frac{\lambda_j^+e^{-\lambda_j^+(t-r)}-\lambda_j^-e^{-\lambda_j^-(t-r)}}
          {\lambda_j^+-\lambda_j^-}\,
        \dd r
         \bigg|^2
\left<x_2,\phi_j\right>^2
          \bigg)^{\frac{1}{2}}.
\end{align}
Further, noticing that
\begin{equation}
0
\leq
\int_s^t
\lambda^{\frac \mu 2} e^{-\lambda(t-r)}\,
\dd r
=\lambda^{\frac{\mu-2}2}(1-e^{-\lambda (t-s)})
\leq
 C(t-s)^{\frac{2-\mu}2},
 \; \text{ for } \; \mu\in[0,2],
 \end{equation}
and recalling (\ref{eq:eqvlent}) again lead us to
\begin{align}
\Big\|
\int_s^t
\mathcal{S}_4(t-r)x_2\,\dd r
      \Big\|
\leq
&
C\bigg(
       \sum_{j=1}^\infty
\bigg|
      \int_s^t\left(e^{-\lambda_j^+(t-r)}
         +
     \lambda_j^-e^{-\lambda_j^-(t-r)}/ \lambda_j
      \right)\,
\dd r
      \bigg|^2
      \left<x_2,\phi_j\right>^2\bigg)^{\frac12}
\nonumber
\\
\leq
&
   2 C \bigg(
             \sum_{j=1}^\infty
     \bigg|
            \int_s^t(\lambda_j^+)^{\frac \mu 2}e^{-\lambda_j^+(t-r)}\,
          \dd r
          \bigg|^2
    (\lambda_j^+)^{-\mu}\left<x_2,\phi_j\right>^2
         \bigg)^{\frac12}
\nonumber\\
&+
    2C \bigg(
    \sum_{j=1}^\infty
      \bigg|
\int_s^t(\lambda_j^-)^{\frac \mu 2}e^{-\lambda_j^-(t-r)}\,\dd r
       \bigg|^2
(\lambda_j^-)^{2-\mu}\lambda_j^{-2}
       \left<x_2,\phi_j\right>^2
          \bigg)^{\frac12}
\nonumber
  \\
\leq&
  C (t-s)^{\frac{2-\mu}{2}}\Big(
       \sum_{j=1}^\infty
                     \lambda_j^{-\mu }\left<x_2,\phi_j\right>^2
       \Big)^{\frac12}
\nonumber \\
  \leq &
C(t-s)^{\frac{2-\mu}2}
       \|x_2\|_{-\mu},
\end{align}
where we also used the fact $\sup_{ j\in \mathbb{N}^+ }(\lambda_j^-)^{2- \mu }<\infty$ for $\mu\in[0,2]$.
This validates Lemma \ref{lem:spatial-regualrity-semigroup}.  $\square$

\begin{lemma}\label{lemma:result-S4}
Let $\mathcal{S}_4(t)$ be a component of the semigroup $\mathcal{S}(t)$ defined in \eqref{eq:definition-Si}.
Then
\begin{align}
\mathbb{J} : =
\bigg\|
     \sum_{l=0}^{m-1}
          \int_{t_l}^{t_{l+1}}
              [\mathcal{S}_4(t_{m}-s)-\mathcal{S}_4(t_{m-l})]
          x_2
    \, \mathrm{d} s
      \bigg\|
      &\leq
     C k\| x_2 \|,
\;
     \forall x_2 \in \dot{H}^0.
\label{eq:result2-S4}
\end{align}
\end{lemma}
{ \it  Proof of Lemma \ref{lemma:result-S4}. }
We apply the same argument  used in the proof of \eqref{lem:integral-regularity2-S4} to arrive at
\begin{align}\label{eq:result-proof-S4}
\mathbb{ J }
\leq
&
 \bigg(
 \sum_{j=1}^\infty
 \bigg|
 \sum_{l=0}^{m-1}
 \int_{t_l}^{t_{l+1}}
     \frac{\lambda_j^+[e^{-\lambda_j^+(t_{m}-s)}-e^{-\lambda_j^+t_{m-l}}]}
                         {\lambda_j^+-\lambda_j^-}\,\dd s
  \bigg|^2
     \left< x_2, \phi_j\right>^2
       \bigg)^{\frac{1}{2}}
       \nonumber
       \\
&+
\bigg(
\sum_{j=1}^\infty
\bigg|
    \sum_{l=0}^{m-1}
      \int_{t_l}^{t_{l+1}}
           \frac{\lambda_j^-[e^{-\lambda_j^-(t_m-s)}-e^{-\lambda_j^-t_{m-l}}]}
                 {\lambda_j^+-\lambda_j^-}
                 \dd s
\bigg|^2
                 \left< x_2, \phi_j\right>^2
      \bigg)^{\frac{1}{2}}.
\end{align}
Note that, for any $\lambda>0$, there exists a constant $C$ such that
\begin{align}
0
\leq
\sum_{l=0}^{m-1}
\int_{t_l}^{t_{l+1}}
[e^{-\lambda(t_m-s)}-e^{-\lambda t_{m-l}}]
\,\dd s
&
=
 \sum_{l=0}^{m-1}
 \int_{t_l}^{t_{l+1}}
\lambda e^{-\lambda(t_m-s)}\lambda^{-1}(1-e^{-\lambda(s-t_l)})
\,\dd s\
\nonumber\\
&\leq
 Ck
 \int_0^{t_m}\lambda e^{-\lambda(t_m-s)}\, \dd s
 \nonumber \\
 & \leq
 Ck.
\end{align}
Hence by inserting this result  into \eqref{eq:result-proof-S4} and  applying \eqref{eq:eqvlent}, one finds that
\begin{align}
\mathbb{ J }
\leq&
 C\bigg(
 \sum_{j=1}^\infty
 \bigg| \sum_{l=0}^{m-1}\int_{t_l}^{t_{l+1}}
      [e^{-\lambda_j^+(t_m-s)}-e^{-\lambda_j^+t_{m-l}}]\dd s
      \bigg|^2\left< x_2,\phi_j\right>^2
      \bigg)^{\frac12}
      \nonumber
      \\
&+
C\bigg(
\sum_{j=1}^\infty
    \bigg|
            \sum_{l=0}^{m-1}
            \int_{t_l}^{t_{l+1}}
              [e^{-\lambda_j^-(t_m-s)}-e^{-\lambda_j^-t_{m-l}}]
            \dd s\bigg|^2
            \left< x_2,\phi_j\right>^2
    \bigg)^{\frac12}
\nonumber
\\ \leq &
Ck \| x_2 \|.
\end{align}
Hence this completes the proof of Lemma \ref{lemma:result-S4}.  $\square$

Armed with the above preparatory results,  we are now ready to prove the main convergence result.

{ \it  Proof of Theorem \ref{them:main-results} }
By recursion, the numerical solution $X^N_m$ can be written as
\begin{align}\label{eq:mild-solution-NM}
X_{m}^N=\mathcal{S}_N(t_{m})X_0^N
+
\sum_{l=0}^{m-1}\int_{t_l}^{t_{l+1}}\mathcal{S}_N(t_{m-l})\mathbf{F}_N(X_l^N)\,\dd s
+
\int_0^{t_m}\mathcal{S}_N(t_m-s)\mathbf{B}_N\,\dd W(s).
\end{align}

Next we deal with the error $u(t_m)-u_m^N$ first.
Owing to \eqref{eq:mild-SPDE}, \eqref{eq:mild-solution-NM},  and the definitions of $\mathcal{S}(t)$ and $\mathcal{S}_N(t)$, one  can easily do the following decompositions,
\begin{align}\label{eq:error-u-um}
u(t_m)-u_m^N
=
&
(I-\mathcal{P}_N)P_1\mathcal{S}(t_m)X_0
+
\sum_{l=0}^{m-1}
\int_{t_l}^{t_{l+1}}
(I-\mathcal{P}_N)\mathcal{S}_2(t_m-s)F(u(s))\,\dd s
\nonumber
\\
&
+
\sum_{l=0}^{m-1}
\int_{t_l}^{t_{l+1}}
[\mathcal{S}_2(t_m-s)-\mathcal{S}_2(t_{m-l})] \mathcal{P}_N F(u(s))
\,\dd s\nonumber
\\
&+
\sum_{l=0}^{m-1}
\int_{t_l}^{t_{l+1}}
\mathcal{S}_{2} ( t_{m-l} ) \mathcal{P}_N [ F(u(s))-F(u(t_l)) ]
\,\dd s
\nonumber\\
&
+
\sum_{l=0}^{m-1}
\int_{t_l}^{t_{l+1}}
\mathcal{S}_{2}(t_{m-l}) \mathcal{P}_N [F(u(t_l))-F(u_l^N)]
\,\dd s\nonumber\\
&
+
\int_0^{t_m}
\mathcal{S}_2(t_m-s)(I-\mathcal{P}_N)\,\dd W(s)
\nonumber\\
:=&
I_1+I_2+I_3+I_4+I_5+I_6.
\end{align}
For the first term $I_1$, note that
 $(u_0, v_0)'\in L^2(\Omega;\dot{H}^{\varrho}\times \dot{H}^{\varrho-2})$ for $\varrho := \min\{1+ \gamma,2\}$, $\gamma\in (-1,2]$. Then  \eqref{lem:spatial-regularity-S1-S2} and \eqref{eq:interpolation-property-PN} help us to get
 \begin{align}
 \|I_1\|_{L^2(\Omega; \dot{H}^0)}
  \leq
  \lambda_{N+1}^{-\frac{\varrho}{2}}
  \|P_1\mathcal{S}(t)X_0\|_{L^2(\Omega;\dot{H}^{\varrho})}
  \leq C\lambda_{N+1}^{-\frac{\varrho}{2}}(\|u_0\|_{L^2(\Omega;\dot{H}^{\varrho})}+
  \|v_0\|_{L^2(\Omega;\dot{H}^{\varrho-2})}) .
 \end{align}
Similarly, using \eqref{eq:interpolation-property-PN},  \eqref{lem:spatial-regularity-F} and
$\eqref{lem:spatial-regularity-S1-S2}$ with $\mu=2$ shows
\begin{align}
\|I_2\|_{L^2(\Omega;\dot{ H}^0)}
&\leq  \lambda_{N+1}^{-1}\int_0^{t_m}
     \|A\mathcal{S}_2( t_m - s )F(u(s))\|_{L^2(\Omega;\dot{H}^0)}\,\dd s\nonumber\\
&\leq
   C\lambda_{N+1}^{-1} \|F(u)\|_{L^\infty([0,T],\;L^2(\Omega;\dot{H}^0))}\nonumber\\
   & \leq
        C\lambda_{N+1}^{-1}(1+\|u_0\|_{L^2(\Omega;\dot{H}^0)}+\|v_0\|_{L^2(\Omega;\dot{H}^{-2})}).
\end{align}
Due to \eqref{lem:spatial-regularity-F} and $\eqref{lemma:temporal-regularity-S1-S2}$ with $\mu=2$,
we obtain the estimate of $I_3:$
\begin{align}
\|I_3\|_{L^2(\Omega; \dot{H}^0)}
&\leq
C \sum_{l=0}^{m-1}
         \int_{t_l}^{t_{l+1}}
             \big\| [\mathcal{S}_2(t_m-s)-\mathcal{S}_2(t_{m-l}) ]P_NF(u (s) ) \big\|_{ L^2(\Omega;\dot{H}^0)}
\,\dd s
\nonumber\\&
\leq
C\sum_{l=0}^{m-1}
\int_{t_l}^{t_{l+1}}
     (s-t_l)
        \,\dd s\,
             \|F( u )\|_{L^\infty([0,T],\;L^2(\Omega;\dot{H}^0))}
        \nonumber\\
&\leq Ck(1+\|u_0\|_{L^2(\Omega;\dot{H}^0)}+\|v_0\|_{L^2(\Omega;\dot{H}^{-2})}).
\end{align}
To treat the term $I_4$, we invoke  a linearization step and the property of $\mathcal{S}_2(t)$ to get
\begin{align}\label{eq:estimation-term-I5}
\|I_4\|_{L^2(\Omega; \dot{H}^0)}
& \leq
\sum_{l=0}^{m-1}\int_{t_l}^{t_{l+1}}
\int_0^1\|\mathcal{S}_{2} (t_{m - l}) P_N F'(\chi(t_l,s, r))[u(s)-u(t_l)]\|_{L^2(\Omega;\dot{H}^0)}\,\dd r\,\dd s
\nonumber \\
& \leq
 C\sum_{l=0}^{m-1}
           \int_{t_l}^{t_{l+1}}
                     \int_0^1
                     \|A^{-\frac{\theta}{2}}F'(\chi(t_l,r,s))[u(s)-u(t_l)]\|_{L^2(\Omega;\dot{H}^0)}
                     \,\dd r\,\dd s,
\end{align}
for $\theta \in [1, 2)$, where for simplicity of presentation we denote
\begin{equation}
\chi (t_l,s, r ):=u(t_l)+r(u(s)-u(t_l)).
\end{equation}
Recall that  $ 1 - \min\{ 1, | \gamma | \} \leq \beta := \min\{1,1+\gamma\} \leq 
\varrho :=\min\{1+ \gamma,2\}$ for $\gamma\in (-1,2]$. 
Hence employing  \eqref{eq:F-assum-condition3}, 
\eqref{eq:thm1-u-spatial2}, \eqref{eq:thm1-u-temporal2}, and H\"older's inequality yields
\begin{align}\label{eq:F(u)-F(uk)}
 \|A^{-\frac{\theta}{2}}F'(\chi(t_l,s,r))[u(s)&-u(t_l)]\|_{L^2(\Omega;\dot{H}^0)}
 \leq
 K (\|\chi(t_l,s,r)\|_{L^4(\Omega;\dot{H}^{\beta})}+1)
        \|u(s)-u(t_l)\|_{L^4(\Omega;\dot{H}^{-\beta})}
        \nonumber\\
     &\leq  C (\|u \|_{L^\infty([0,T],\;L^4(\Omega;\dot{H}^{\varrho}))}+1)
        \|u(s)-u(t_l)\|_{L^4(\Omega;\dot{H}^{-\delta})}
        \nonumber\\
& \leq
  C k^{\min\{1+\gamma,1\}}
      \big(1 + \|u_0\|^2_{L^4(\Omega;\dot{H}^{\varrho})}
      +
      \|v_0\|^2_{L^4(\Omega;\dot{H}^{\varrho-2})} \big).
\end{align}
Plugging this into \eqref{eq:estimation-term-I5} yields
\begin{align}\label{eq:estimation-I5}
\|I_4\|_{L^2(\Omega; \dot{H}^0)}
\leq
 C k^{\min\{1+\gamma,1\}}
            (1+ \|u_0\|^2_{L^4(\Omega;\dot{H}^{\varrho})}
            +
            \|v_0\|^2_{L^4(\Omega;\dot{H}^{\varrho-2})}).
\end{align}
Next, with the aid of \eqref{lem:spatial-regularity-S1-S2} and the Lipschitz condition \eqref{eq:F-assum-condition1}, we derive for  the term $I_5$ that
\begin{align}
\|I_5\|_{L^2(\Omega; \dot{H}^0)}
&\leq
\sum_{l=0}^{m-1}
\int_{t_l}^{t_{l+1}}
\|\mathcal{S}_{2 } ( t_{m-l} ) \mathcal{P}_N \|_{\mathcal{L}(\dot{H}^0)}\|F(u(t_l))-F(u_l^N)\|_{L^2(\Omega;\dot{H}^0)}
\,\dd s
\nonumber\\
&\leq
C k \sum_{l=0}^{m-1} \|u(t_l)-u_l^N\|_{L^2(\Omega;\dot{H}^0)}.
\end{align}
For the last term $I_6$, we use \eqref{lem:spatial-regularity-S1-S2},
the It\^{o} isometry and  the property of $P_N$ to show
\begin{align}
\|I_6\|_{L^2(\Omega; \dot{H}^0)}
&=
\Big(\int_0^{t_m}
\|
   \mathcal{S}_2(t_m-s)(I-\mathcal{P}_N)Q^{\frac12}
         \|_{\mathrm{HS}}^2
         \,\dd s\Big)^{\frac12}
\nonumber\\&
\leq
C\sqrt{T}
    \|A^{-1}(I-\mathcal{P}_N)Q^{\frac12}\|_{\mathrm{HS}}\nonumber\\
 &\leq
     C\lambda_{N+1}^{-\frac\varrho2}
            \|A^{\frac{\varrho-2}2}Q^{\frac12}\|_{\mathrm{HS}}
   \nonumber\\          
  & \leq
     C\lambda_{N+1}^{-\frac{\varrho}2}
            \|A^{\frac{\gamma-1}2}Q^{\frac12}\|_{\mathrm{HS}}.
\end{align}
Eventually, inserting the above estimates into \eqref{eq:error-u-um}
and applying the discrete Gronwall's inequality imply that,  for $\varrho = \min\{1+ \gamma,2\}$,
\begin{align}\label{eq:Error-u-um}
\|u(t_m)-u_m^N\|_{L^2(\Omega;H)}\leq C(\lambda_{N+1}^{-\frac{\varrho}{2}}+k^{\min\{1+\gamma,1\}})
      \big(1 + \|u_0\|^2_{L^4(\Omega;\dot{H}^{\varrho})}
      +
      \|v_0\|^2_{L^4(\Omega;\dot{H}^{\varrho-2})} \big),
\end{align}
which verifies \eqref{them:error-estimates-u-um}.
Next we start to treat the error $\|v(t_m)-v_m^N\|_{L^2(\Omega;\dot{H}^0)}$:
\begin{align}
v(t_m) - v_m^N
& =
(I-\mathcal{P}_N)P_2\mathcal{S}(t_m)
X_0
+
\int_0^{t_m}
\mathcal{S}_4(t_m-s)(I-\mathcal{P}_N)
\,\dd W(s)
\nonumber
\\
& \quad
+
\sum_{l=0}^{m-1}
\int_{t_l}^{t_{l+1}}
\big[
\mathcal{S}_4(t_m-s)F(u(s)) - \mathcal{S}_{4,N}(t_{m-l})\mathcal{P}_NF(u_l^N)
\big]
\,\dd s
\nonumber\\
& :=
J_1 + J_2 + J_3.
\end{align}
Analogously to the estimate of $I_1$, we use \eqref{lem:spatial-regularity-S3-S4} with $\mu = 0$ to bound the term $J_1$:
\begin{align}\label{eq:estiamtion-L1}
\|J_1\|_{L^2(\Omega; \dot{H}^0)}
\leq
C\lambda_{N+1}^{-\frac\gamma2}
\big(
\|u_0\|_{L^2(\Omega;\dot{H}^\gamma)}+\|v_0\|_{L^2(\Omega;\dot{H}^\gamma)}
\big)
\quad \text{ for } \gamma\in [0,2].
\end{align}
Owing to \eqref{lem:integral-regularity1-S4} with $\nu = 0$,
Assumption \ref{assum:wiener-process-condition} and the It\^{o} isometry, we acquire that,  for $\gamma\in [0,2]$,
\begin{align}\label{eq:estimation-L3}
\|J_2\|_{L^2(\Omega; \dot{H}^0)}
&=
\Big(
\int_0^{t_m}
\|\mathcal{S}_4(t_m-s)(I-\mathcal{P}_N)Q^{\frac12}\|^2_{\mathrm{HS}}
\,\dd s
\Big)^{\frac12}
\nonumber\\
&
\leq
 C\|A^{-\frac12}(I-\mathcal{P}_N)Q^{\frac12}\|_{\mathrm{HS}}
\nonumber \\
& \leq
 C\lambda_{N+1}^{-\frac{\gamma}{2}}
    \|A^{\frac{\gamma-1}2}Q^{\frac12}\|_{\mathrm{HS}}.
\end{align}
Now it remains to estimate $J_3$. To this aim, we furthermore decompose $J_3$ into three terms
\begin{align}
\|J_3 \|_{L^2(\Omega; \dot{H}^0)}
\leq&
\Big\|
      \sum_{l=0}^{m-1}\int_{t_l}^{t_{l+1}}(I-\mathcal{P}_N)\mathcal{S}_4(t_m-s)
         F(u(s))\,\dd s
       \Big\|_{L^2(\Omega;\dot{H}^0)}
\nonumber\\
&+
\Big\|
      \sum_{l=0}^{m-1}
         \int_{t_l}^{t_{l+1}}
            \mathcal{P}_N[\mathcal{S}_4(t_m-s)-\mathcal{S}_4(t_{m-l})]
                 F(u(s))\, \dd s
        \Big\|_{L^2(\Omega;\dot{H}^0)}
\nonumber\\
&+
\Big\|
\sum_{l=0}^{m-1}
    \int_{t_l}^{t_{l+1}}
          \mathcal{P}_N\mathcal{S}_4(t_{m-l})[F(u(s))-F(u_l^N)]\,\dd s
        \Big\|_{L^2(\Omega;\dot{H}^0)}
\nonumber\\
:= &
J_3^1+ J_3^2 + J_3^3.
\end{align}
Below, we will estimate them separately. Concerning $J_3^1$,  we employ \eqref{lem:temporal-regularity-F}, \eqref{lem:spatial-regularity-S3-S4} with $\mu = 2$, \eqref{lem:integral-regularity2-S4} with $\mu = 2$
and \eqref{eq:interpolation-property-PN}  to arrive at, for $\varrho= \min\{1 + \gamma,2\}$, $\gamma\in[0,2]$,
\begin{align}\label{eq:estimaion-L21}
J_3^1
\leq&
 \lambda_{N+1}^{-1}  \Big[\int_0^{t_m} \left\| A\mathcal{S}_4(t_m-s)
      [F(u(s))-F(u(t_m))]\right\|_{L^2(\Omega;\dot{H}^0)}   \,\dd s
\nonumber\\&
+
\Big\|
     \int_0^{t_m}A\mathcal{S}_4(t_m-s)
          F(u(t_m))\,\dd s
      \Big\|_{L^2(\Omega;\dot{H}^0)}\Big]
\nonumber\\
\leq&
  C\lambda_{N+1}^{-1}\Big[ \int_0^{t_m}
       (t_m-s)^{-1}
                   \|F(u(s))-F(u(t_m))\|_{L^2(\Omega;\dot{H}^0)}\,\dd s
+
             \|F(u(t_m))\|_{L^2(\Omega;\dot{H}^0)}\Big]
\nonumber\\
\leq&
 C \lambda_{N+1}^{-1}\Big[\int_0^{t_m}
         (t_m-s)^{-1}
                        (t_m-s)^{\frac{\varrho} {2} }\,\dd s
 +
 1\Big]
    (1+\|u_0\|_{L^2(\Omega; \dot{H}^\varrho)}
    +
    \|v_0\|_{L^2(\Omega; \dot{H}^{\varrho-2})})\nonumber\\
 \leq&
     C\lambda_{N+1}^{-1}
     (1+\|u_0\|_{L^2(\Omega; \dot{H}^\varrho)}
     +
     \|v_0\|_{L^2(\Omega; \dot{H}^{\varrho-2})}).
\end{align}
Further,  we  utilize \eqref{lem:temporal-regularity-F}, \eqref{lemma:temporal-regularity-S4} with $\mu=2$, $\nu=0$ and
 \eqref{eq:result2-S4} with $x_2 =F(u(t_m))$ to achieve
\begin{align}\label{eq:estimation-L22}
J_3^2
\leq
&
\sum_{l=0}^{m-1}
      \int_{t_l}^{t_{l+1}}
     \left\|
     [\mathcal{S}_4(t_m-s)-\mathcal{S}_4(t_{m-l})]
\left[
F(u(s))-F(u(t_m))
      \right]
              \right\|_{L^2(\Omega;\dot{H}^0)}\,\dd s
              \nonumber
      \\
&+
\Big\|
      \sum_{l=0}^{m-1}
           \int_{t_l}^{t_{l+1}}
                 [\mathcal{S}_4(t_m-s)-\mathcal{S}_4(t_{m-l})]
            F(u(t_m))
     \,\dd s
          \Big\|_{L^2(\Omega;\dot{H}^0)}
     \nonumber\\
\leq&
 Ck\Big[\sum_{l=0}^{m-1}
     \int_{t_l}^{t_{l+1}}
       (t_m-s)^{-1}
                \|F(u(s))-F(u(t_m))\|_{L^2(\Omega;\dot{H}^0)}
    \, \dd s
 +
     \|F(u(t_m))\|_{L^2(\Omega;\dot{H}^0)}\Big]
\nonumber\\
\leq &
Ck\Big[
\int_0^{t_m}(t_m-s)^{-1}(t_m-s)^{\frac{\varrho} {2} }\,\dd s
   +1
   \Big]
   (1+\|u_0\|_{L^2(\Omega; \dot{H}^\varrho)}
   +
   \|v_0\|_{L^2(\Omega; \dot{H}^{\varrho-2})})
\nonumber\\
   \leq &
   Ck(1+\|u_0\|_{L^2(\Omega; \dot{H}^\varrho)}
   +
   \|v_0\|_{L^2(\Omega; \dot{H}^{\varrho-2})}).
\end{align}
At last, similarly to the proof  of \eqref{eq:estimation-I5}, one can employ \eqref{lem:spatial-regularity-S3-S4}
with $\mu = \theta$,  \eqref{eq:Error-u-um}, \eqref{eq:F-assum-condition3} and H\"older's inequality
to show that, for $\theta\in [1,2)$ and $\gamma\in [0,2]$,
\begin{align}\label{eq:estimation-L23}
J_3^3
\leq
&
\sum_{l=0}^{m-1}
\int_{t_l}^{t_{l+1}}
\Big[\| \mathcal{S}_4(t_{m-l})[F(u(s))-F(u(t_l))]\|_{L^2(\Omega;\dot{H}^0)}
 +
                 \| \mathcal{S}_4(t_{m-l})[F(u(t_l))-F(u_l^N)] \|_{L^2(\Omega;\dot{H}^0)} \Big]
\,\dd s
\nonumber\\
\leq
&
C\sum_{l=0}^{m-1}
\Big[
      \int_{t_l}^{t_{l+1}}
       t_{m-l}^{-\frac\theta2}
                \|A^{-\frac\theta2} [F(u(s))-F(u(t_l))] \|_{L^2(\Omega;\dot{H}^0)}
\,\dd s
+
      k\|u(t_l)-u_l^N\|_{L^2(\Omega;\dot{H}^0)}
     \Big]
\nonumber\\
\leq&
 C
 k\sum_{l=0}^{m-1}
            \Big[k^{\min\{1+\gamma,1\}}
                t_{m-l}^{-\frac\theta 2 }
 +
  \lambda_{N+1}^{-\frac{1+\min\{\gamma, 1\}}2  }+k^{\min\{1+\gamma,1\}})
      \Big]
     (1 + \|u_0\|^2_{L^4( \Omega;\dot{H}^{\varrho} )}
     +
     \|v_0\|^2_{ L^4( \Omega;\dot{H}^{\varrho-2} ) })
     \nonumber
     \\
\leq&
 C(k+\lambda_{N+1}^{- \frac{ 1+\min\{\gamma, 1\} }2 })
(1+\|u_0\|^2_{ L^4( \Omega;\dot{H}^{\varrho} ) }
+
\|v_0\|^2_{ L^4( \Omega;\dot{H}^{ \varrho-2 } ) } ).
\end{align}
Gathering \eqref{eq:estimaion-L21}, \eqref{eq:estimation-L22} and \eqref{eq:estimation-L23} together gives
\begin{align}
\|J_3\|_{L^2(\Omega;\dot{H}^0)}
\leq
   C(k+\lambda_{N+1}^{-\frac{\gamma}{2}})
    (1+\|u_0\|_{L^4(\Omega;\dot{H}^{1+\min\{\gamma,1\}})}+\|v_0\|_{L^4(\Omega;\dot{H}^{\gamma})})^2,
\end{align}
which together with \eqref{eq:estiamtion-L1} and \eqref{eq:estimation-L3} yields
\begin{align}
\|v(t_m)-v_m^N\|_{L^2(\Omega;\dot{H}^0)}\leq C(k+\lambda_{N+1}^{-\frac{\gamma}{2}})
    (1+\|u_0\|_{L^4(\Omega;\dot{H}^{1+\min\{\gamma,1\}})}+\|v_0\|_{L^4(\Omega;\dot{H}^{\gamma})})^2.
\end{align}
The proof of Theorem \ref{them:main-results} is finally complete. $\square$

\section{Numerical experiments}
\label{sec:numer-exp}



In this section, we test the previous theoretical findings by doing some numerical experiments. For simplicity,
we take a stochastic strongly damped wave equation in one  dimension as follows
\begin{align}\label{eq:one-dimension-steq}
\left\{\begin{array}{ll}
u_{tt}
=
\Delta u
+
\Delta u_t + \frac {1-u}{1+u^2}+\dot{W}(t),
&
t\in(0,1],
\;
x\in(0,1),
\\
u(0,x)=\frac{\partial u}{\partial t}(0,x)=0,
&
x\in (0,1),
\\
u(t,0)=u(t,1)=0,
&
t\in(0,1].
\end{array}\right.
\end{align}
Our task is to simulate the approximation errors at a fixed time $T=1$
for the space-time white noise case  ($Q=I$) and the trace-class noise case  ($Q=A^{-0.5005}$).
The numerical errors in the mean-square sense are achieved by the Monte-Carlo approach over 100 samples.
Since no exact  solution is available,  we identify it  as a numerical one with
fine mesh step-sizes by choosing small $k_{\text{exact} }$ and large $N_{ \text{exact} }$.
\begin{figure}[!ht]
\centering
      \includegraphics[width=3in,height=2.8in] {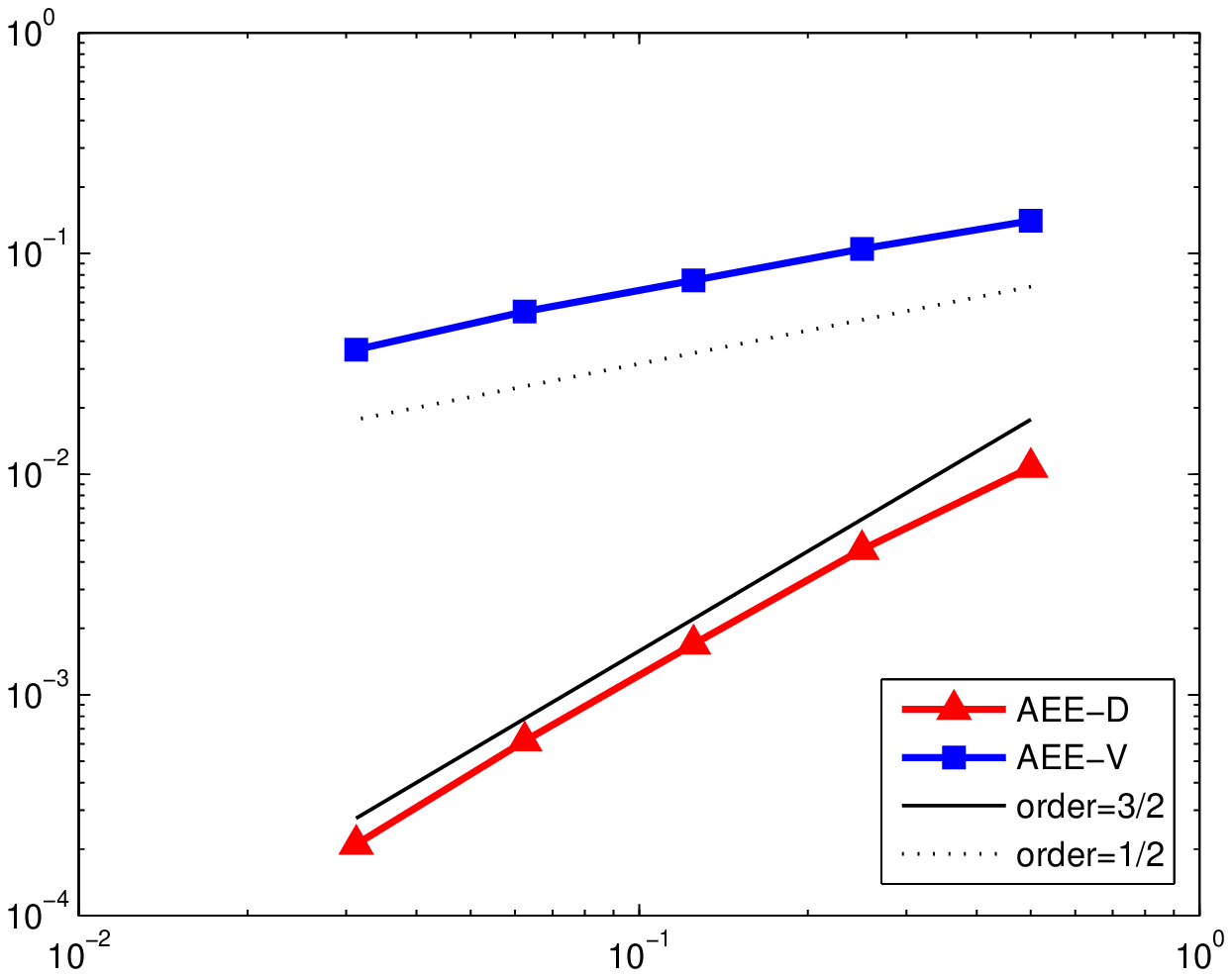}
       \includegraphics[width=3in,height=2.8in] {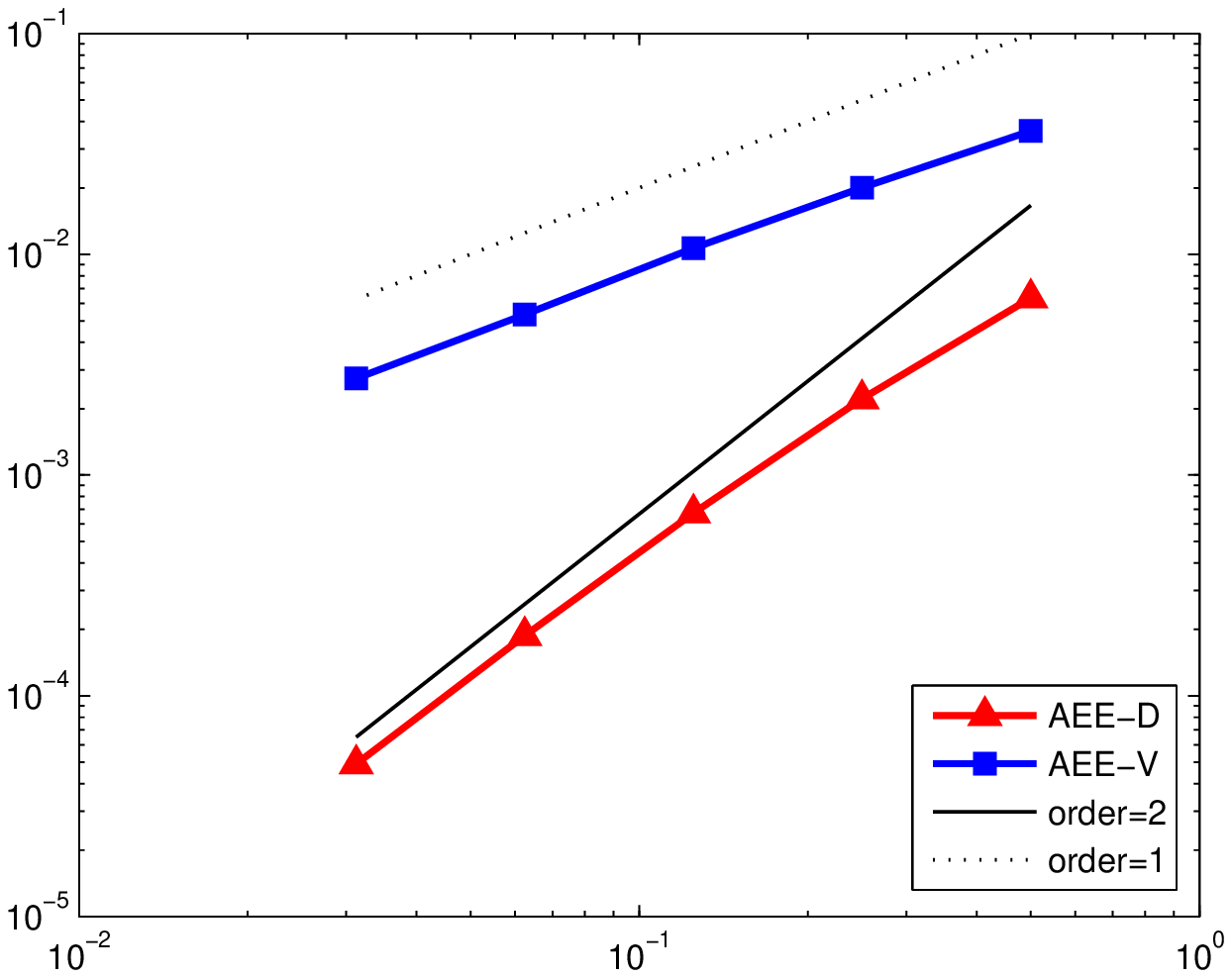}
 \caption{Convergence rates for the spatial discretization (Left: $Q=I$; right: $Q=A^{-0.5005})$.}
\label{FG:spatial-error}
\end{figure}

At first, we perform the spatial discretization of the test equation with two noise cases mentioned earlier.
Fig.\ref{FG:spatial-error} depicts the spatial numerical errors of the accelerated exponential Euler (AEE) scheme \eqref{eq:full-scheme-form} with  $N=2^j, j = 1, 2,...,5$ and  time step-size $k_{\text{exact} }=2^{-14}$.
Evidently, the resulting errors for both the displacement (AEE-D) and the velocity (AEE-V)
decrease with slopes of orders as expected (see  Theorem \ref{them:main-results}).
Here note that the "exact" solution is computed by using $N_{\text{exact} }=2^8$.

\begin{figure}[!ht]
\centering
      \includegraphics[width=3in,height=3in] {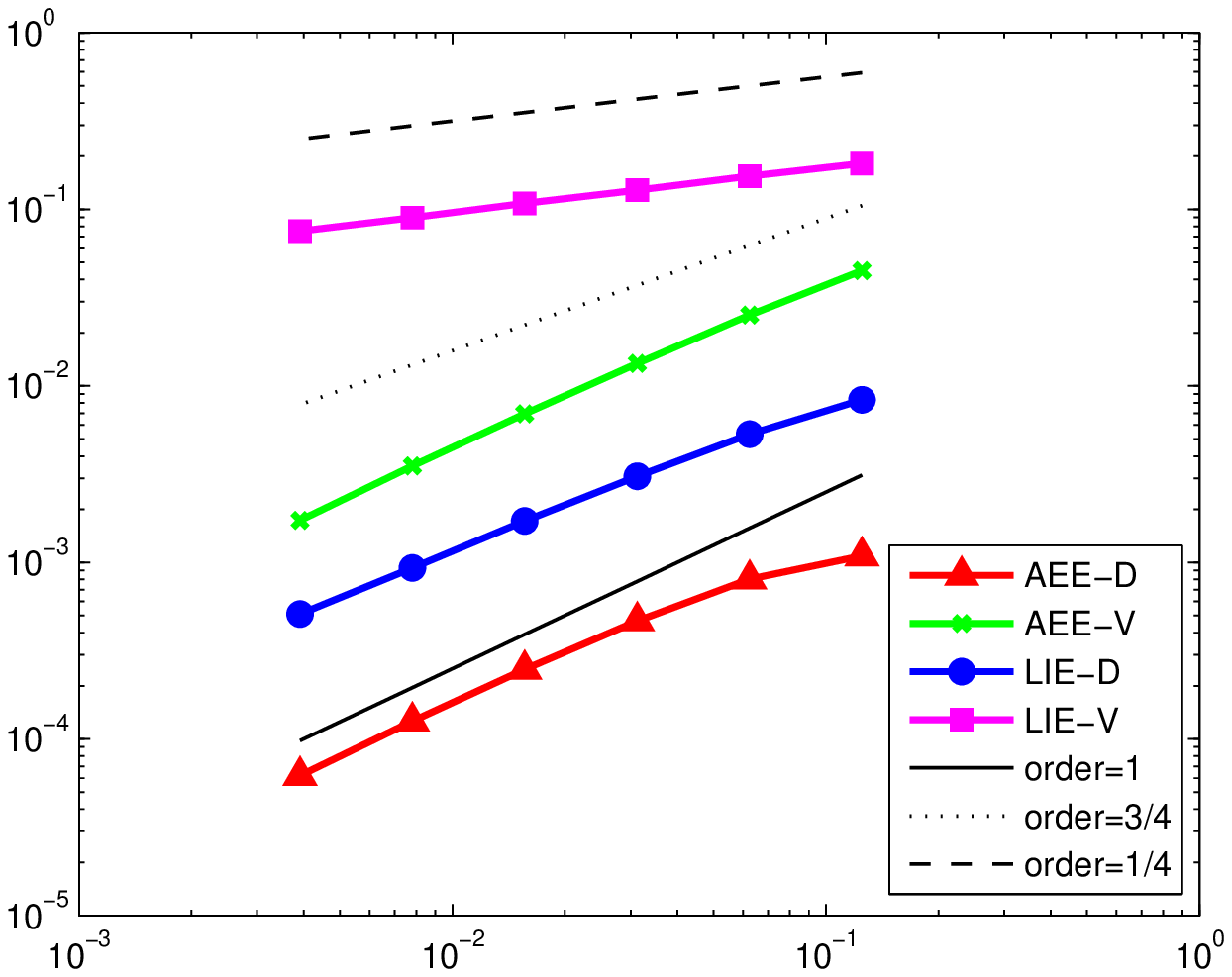}
       \includegraphics[width=3in,height=3in] {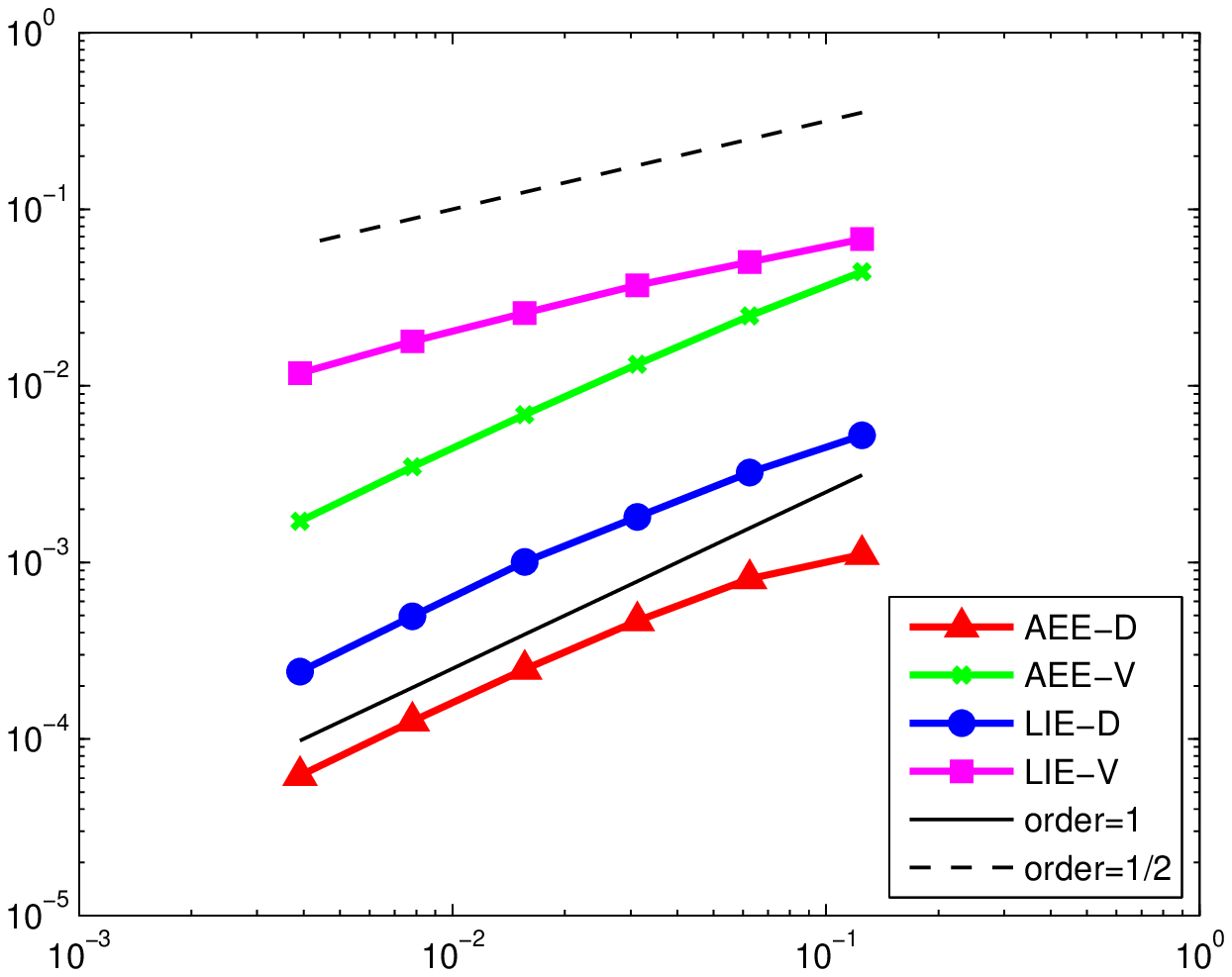}
 \caption{Convergence rates for the temporal discretization (Left: $Q=I$; right: $Q=A^{-0.5005}$). }
  \label{FG:time-error}
\end{figure}

In what follows, we examine the temporal discretization errors
using different time stepsizes $k = 2^{-i}, i=3,4,...,8$. The "exact" solution is
approximated by AEE scheme \eqref{eq:full-scheme-form} with small time step-size $k_{\text{exact} }=2^{-12}$.
For comparison, we present in Fig.\ref{FG:time-error} the convergence errors in time caused
by the AEE scheme \eqref{eq:full-scheme-form}  and the linear implicit Euler (LIE) scheme  \cite{qi2015error}.
From Fig.\ref{FG:time-error}, one can observe that, the errors of \eqref{eq:full-scheme-form} for the displacement (AEE-D) and  for the velocity (AEE-V) both decrease with order $1$ for the space-time 
white noise and the trace-class noise, which is consistent with assertions in 
Theorem \ref{them:main-results}. The corresponding errors for the LIE method 
(LIE-D and LIE-V), however, exhibit much worse  performance.  For example, the errors of LIE-D and LIE-V 
only show rates of order $\tfrac34$ and  order $\tfrac14$, respectively, in the space-time white noise case.
It turns out that the newly proposed scheme admits significantly improved convergence rates in time.

\bibliography{bibfile}

\bibliographystyle{abbrv}
 \end{document}